\documentclass[12pt]{article}
\usepackage{graphics,natbib}
\usepackage[scanall]{psfrag}
\usepackage[colorlinks=true,citecolor=blue]{hyperref}
\usepackage{graphics}
\usepackage{natbib}
\bibpunct{(}{)}{;}{x}{,}{,}

\usepackage{amsmath}
\usepackage{amssymb}
\usepackage{fancybox}
\usepackage{fancyheadings}

\newcommand{\be}{\begin{equation}}
\newcommand{\ee}{\end{equation}}
\renewcommand{\k}{\kappa}

\newcommand{\M}{\mathcal{M}}
\renewcommand{\P}{\mathcal{P}}
\newcommand{\F}{\mathcal{F}}
\renewcommand{\H}{\mathcal{H}}
\newcommand{\G}{\mathcal{G}}

\newcommand{\Q}{\mathcal{Q}}
\renewcommand{\S}{\mathcal{S}}
\newcommand{\R}{\mathcal{R}}
\newcommand{\mT}{\mathcal{T}}
\newcommand{\X}{\mathcal{X}}
\newcommand{\mE}{\mathbb{E}}

\newcommand{\mC}{\mathbb{C}}
\newcommand{\mP}{\mathbb{P}}
\newcommand{\mV}{\mathbb{V}}

\newcommand{\mz}{\mathfrak{Z}}
\newcommand{\T}{\mathfrak{T}}

\newcommand{\V}{V}

\newcommand{\bx}{\mathbf{x}}
\newcommand{\bU}{\mathbf{U}}
\newcommand{\bQ}{\mathbf{Q}}
\newcommand{\bJ}{\mathbf{J}}

\newcommand{\mbb}{\boldsymbol{\beta}}
\newcommand{\mbt}{\boldsymbol{\theta}}
\newcommand{\mbtu}{\boldsymbol{\tau}}
\newcommand{\mbl}{\boldsymbol{\lambda}}

\newcommand{\mbL}{\boldsymbol{\Lambda}}
\newcommand{\mbx}{\boldsymbol{\xi}}
\newcommand{\mbu}{\boldsymbol{\mz}}
\newcommand{\mbT}{\boldsymbol{\Theta}}

\newcommand{\matr}[1]{\mathbf{#1}}
\newcommand{\ip}[2]{\left <{#1},{#2}\right >}
\newcommand{\paren}[1]{\left({#1}\right)}
\newcommand{\cov}{\mbox{cov}}
\newcommand{\cor}{\mbox{cor}}

\newcommand{\sgn}{\mbox{sgn}}
\newcommand{\supp}{\mbox{supp}}
\newcommand{\toinlaw}{\rightsquigarrow}

\newcommand{\toinas}{\stackrel{a.s.}{\longrightarrow}}

\newcommand{\we}{\widehat{\eta}}
\newcommand{\wmbt}{\widehat{\mbt} }
\newcommand{\wQ}{\widehat{\Q}}

\newcommand{\supt}{\sup_{\mbt \, \in \, \mbT}}
\newcommand{\inft}{\inf_{\mbt \, \in \, \mbT}}

\newcounter{assumption}
\renewcommand{\theassumption}{A\arabic{assumption}}
\newcommand{\assumption}{
  \refstepcounter{assumption}
  \noindent\textbf{\theassumption:} }

\newcounter{remark}
\renewcommand{\theremark}{\arabic{remark}}
\newcommand{\remark}{
  \refstepcounter{remark}
  \noindent{\em Remark \theremark:} }

\newtheorem{theorem}{Theorem}
\newtheorem{lemma}{Lemma}

\begin{document}

\begin{center}
{\large \bf Inference in Perturbation Models, Finite Mixtures and Scan
Statistics: The Volume-of-Tube Formula} \\

Ramani S.\ Pilla\footnote{Research supported in part by the
  National Science Foundation (NSF) grant DMS 02-39053 and
  Office of Naval Research (ONR) grants N00014-02-1-0316 and
  N00014-04-1-0481.} and Catherine Loader\footnote{Research supported
  in part by the NSF grant DMS 03-06202 and ONR grant
  N00014-04-1-0481.}  \\
Department of Statistics, Case Western Reserve University, Cleveland,
OH 44106 \\ 
pilla@case.edu \hspace{0.05in} catherine@case.edu \\
{\bf Abstract}
\end{center}
This research creates a general class of {\em perturbation models}
which are described by an underlying {\em null} model that accounts
for most of the structure in data and a perturbation that accounts for
possible small localized departures. The perturbation models encompass
finite mixture models and spatial scan process.  In this article, (1)
we propose a new test statistic to detect the presence of
perturbation, including the case where the null model contains a set
of nuisance parameters, and show that it is equivalent to the
likelihood ratio test; (2) we establish that the asymptotic
distribution of the test statistic is equivalent to the supremum of a
Gaussian random field over a high-dimensional manifold (e.g., curve,
surface etc.) with boundaries and singularities; (3) we derive a
technique for approximating the quantiles of the test statistic using
the Hotelling-Weyl-Naiman {\em volume-of-tube formula}; and (4) we
solve the long-pending problem of testing for the order of a mixture
model; in particular, derive the asymptotic null distribution for a
general family of mixture models including the multivariate
mixtures. The inferential theory developed in this article is
applicable for a class of non-regular statistical problems involving
loss of identifiability or when some of the parameters are on the
boundary of the parametric space.

\noindent{\em Keywords:} Gaussian random field, Likelihood ratio test
statistic, Multivariate Mixture Models, Nonparametric maximum
likelihood estimator, Nuisance parameters, Score process,
Volume-of-tube formula.

\setcounter{equation}{0} \def\theequation{\thesection.\arabic{equation}}

\section{Introduction and Motivation}
\label{sec-intro}

A fundamental and yet a very challenging problem in finite mixtures is
determining the {\em order of a mixture model} or {\em mixture
complexity}. This problem has been under intense investigation for
over thirty years \citep{wolfe:71,roeder:90,lindsay:95} with no
practically feasible solution for a general class of mixture
families. Establishing a valid large-sample theoretical framework
along with a practically feasible machinery for testing the order of a
mixture model formed from a broad class of densities remains an open
problem and is the focus of this research. It has long been noted that
testing for the number of mixture components is a non-regular problem 
(a) due to loss of identifiability of the null distribution (i.e., the
parameters representing the null distribution are not unique) and (b)
since the parameters under the null hypothesis are on the boundary of
the parameter space, instead of its interior. Consequently, the
likelihood ratio test (LRT) statistic does not have the standard
asymptotic null distribution of chi-squared 
\citep{chernoff:54,ghosh:85,hart:85,bickel:93}. As noted by several  
authors, the asymptotic null distribution of the LRT statistic is
highly complex and very difficult to simulate from in practice.

The main thrust of this research is to create a fundamental class of models
referred to as {\em perturbation models} and derive large-sample
theory to detect the presence of perturbation. These models play
an instrumental role in the development of inferential theory for a
class of important problems such as (1) testing for the order of a
mixture model formed from smooth families of densities, including the
multivariate case; (2) searching for an unusual activity or region in
the context of spatial scan process; and (3) detecting a signal in
the presence of noisy backgrounds \citep{pilla:05}. The resulting
theory has broad applications in astronomy, astrophysics, biology,
medicine, particle physics and datamining, to name a few.

\subsection{Perturbation Models}  

Let $\P = \{p(x; \eta, \mbl, \mbt)\!: \mbl \in \mbL, \mbt \in \mbT 
\subset \R^d\}$ be a family of probability density functions. Assume
that ${\bf X} = (X_1, \ldots, X_n)$ is an independently and
identically distributed (i.i.d.) random sample from
\begin{equation}
  \label{eq:pbm}
  p(x; \eta, \mbl, \mbt) := (1 - \eta) \, f(x; \mbl) + \eta \, \psi(x;
  \mbt),   
\end{equation} 
where \(f(\cdot; \mbl)\) is a {\em null density} for an unknown
parameter vector $\mbl \in \mbL$, \(\psi(\cdot; \mbt)\) is a {\em
  perturbation density} with an unknown {\em nuisance parameter}
vector \(\mbt \in \mbT \subset \R^d\), both defined on a sample space
$\X \subset \R^s$ and \(\eta \in [0, 1]\) is the size of the
perturbation. In the context of finite mixture models, the null model
represents a mixture with \(m\) component densities and the
perturbation model represents additional component densities. In the
spatial scan process scenario, the null density accounts for the
background or noise whereas the perturbation searches for an unusual
activity.

The central idea is to introduce a {\em perturbation parameter} $\eta$
which creates a departure from the null model. There are two primary
goals: (1) Estimation of the parameters in the perturbation model and
(2) testing the hypothesis 
\begin{eqnarray}
  \label{eq:bashyp}
  \H_0\!: \eta = 0 \quad \mbox{against} \quad  \H_1\!: \eta > 0.
\end{eqnarray}
Under $\H_0$, $p(\cdot; \eta, \mbl, \mbt) = f(\cdot; \mbl)$ and the
null model entirely describes the data. However, under $\H_1$, the
term $\eta \, \psi(\cdot; \mbt)$ represents a departure from the null
model.

The perturbation model falls into a class of problems studied by
\cite{davies:77,davies:87} in which a vector of {\em nuisance
  parameters} (in our case \(\mbt\)) appears only under the
alternative hypothesis and standard asymptotic theory for the LRT
breaks down. In particular, the asymptotic behavior of the LRT for the
testing problem (\ref{eq:bashyp}) is very difficult to characterize
due to the difficulties with the geometry of the parameter space
(scenarios (a) and (b) discussed earlier). It is worth noting that
these same set of problems occur in the context of testing for
homogeneity in finite mixture models. The inferential theory developed
in this article requires only mild smoothness conditions on the family
of densities while being generic and applicable much more widely. The
two most important and distinct statistical problems motivating this
work are finite mixture models \citep{lindsay:95} and spatial scan
analysis \citep{glaz:01}. 

\subsection{Inference in Mixture Models} 
\label{sec-mix}

Let $\F = \{\psi(x; \mbt): \mbt \in \mbT \subset \R^d\}$ be a family
of probability densities with respect to a \(\sigma\)-finite
dominating measure $\mu$ for an $s$-dimensional random vector $x \in
\X \subset \R^s$ and let $\G$ be the space of all probability measures
on $\mbT$ with the $\sigma$-field generated by its Borel
subsets. Assume that the component density $\psi(\cdot; \mbt)$ is bounded
in $\mbt \in \mbT$. 

Suppose that given \(\mbt\), a random variable \(X\) has a density
\(\psi(x; \mbt)\) and that $\mbt$ follows a distribution $\Q$,
referred to as {\em mixing distribution}. For a given $\Q \in \G$,
assume that the sample arises from the marginal density ${\rm g}(x;
\Q) := \int_{\mbt} \psi(x; \mbt) \; d \, \Q(\mbt)$ for $x \in \X
\subset \R^{s}$ referred to as a {\em mixture density} with a
corresponding {\em mixing measure} $\Q$.  In the case of a discrete
and finitely supported mixing measure, the {\em mixing distribution}
can be expressed as $\Q_m = \sum_{j = 1}^m \beta_j \,
\varepsilon(\mbt_j)$, where $\varepsilon(\cdot)$ is a point mass
function and $\mbt_{1}, \ldots, \mbt_{m}$ are distinct {\em support
  point vectors} with a corresponding vector of {\em mixing weights}
$\boldsymbol{\beta} := (\beta_{1}, \ldots, \beta_{m})^T$ such that
$\boldsymbol{\beta}$ belongs to the interior of the unit simplex
$\{\boldsymbol{\beta}\!: \sum_{j = 1}^m \beta_j = 1, \beta_j \geq 0, j
= 1, \ldots, m\}$. Therefore, mixture density can be expressed as
${\rm g}(x; \Q_m) = \sum_{j = 1}^m \beta_{j} \, \psi(x; \mbt_j)$,
where the number of support points $m$ becomes the order of the
mixture model or mixture complexity. The probability distribution
$\Q_m$ that maximizes the loglikelihood $l(\Q_m) = \sum_{i = 1}^n \log
\, [{\rm g}(x_i; \Q_m)]$ is the {\em nonparametric maximum likelihood
  estimator} (NPMLE) of $\Q_m$ \citep{lindsay:95}. 

A long-pending and very challenging problem is determining the order
$m$ of the mixture model. In the perturbation model framework, if $f(;
\mbl)$ represents the \(m\)-component mixture density ${\rm g}(\cdot;
\Q_m)$, then \(\psi(\cdot; \mbt_{m + 1})\) represents the $(m + 1)$st
component density. Therefore, inferential theory for perturbation
models provides the machinery for testing the order of a mixture
model. If $m$ is fixed, the loglikelihood has multiple local maxima
and the LRT has an unknown limiting distribution. In the case of
normal mean mixtures and under severe identifiability conditions,
\cite{ghosh:85} derived the asymptotic null distribution of the LRT as
\begin{eqnarray}
 \label{eq:lrt}
 \supt \, [Z(\mbt)]^2 \; {\bf 1} \, [Z(\mbt) \geq 0], \end{eqnarray} where
 $Z(\mbt)$ is a zero mean Gaussian process indexed by a set $\mbt$
 with a specified covariance function and ${\bf 1}[\cdot]$ is the
 indicator function. When the support set of certain parameters in the
 model is unbounded (e.g., in normal and gamma mixtures), the LRT
 statistic can diverge to infinity as $n \to \infty$ instead of having
 a limiting distribution \citep{hart:85,liu:03}. This divergence of
 the LRT poses major difficulties in characterizing the distribution
 of the LRT and in obtaining reliable simulation results for the null
 distribution \citep{lindsay:95}. For testing in multinomial mixture
 models, \citet{lindsay:95} derived approximation to the asymptotic
 distribution of the LRT based on the Hotelling-Weyl
 \citep{hotelling:39,weyl:39} volume-of-tube formula.

Existing theoretical results have been obtained only for some
 special cases and many researchers have considered simulation and
 resampling based approaches to approximate the asymptotic null
 distribution of the LRT for simple models; see \cite{lindsay:95} and
 \cite{mclach:00} for detailed discussion and other
 references. \citet{dacunha:99} proposed a general theory for the
 asymptotic null distribution of the LRT in testing for $\H_0\!: m =
 p$ mixtures against $\H_1\!: m = q$ mixtures, where $q > p$ using a
 locally conic parameterization. Under certain stringent conditions,
 they showed that the asymptotic null distribution of the LRT
 statistic has a form similar to (\ref{eq:lrt}); however, tail
 probability calculations required for calibrating the LRT statistic
 are not derived. Unfortunately, analytic derivations of the
 distribution of supremum of the Gaussian process are difficult
 problems. Most importantly, the issue of ``singularities of the
 process'' (as described in Section \ref{sec-singular}) is of
 fundamental importance in the context of mixture testing problem and
 it has not been addressed in the existing literature, including by 
 \citet{dacunha:99}.

 The perturbation theory developed in this article, provides an
 elegant and flexible machinery for approximating the quantiles of the
 test statistic for the following class of fundamental problems: (1)
 testing problems in which the true parameter is on the boundary of
 the hypotheses regions; (2) testing $\H_0\!: \mbox{$m$-component
   mixture}$ against $\H_1\!: \mbox{$(m + q)$-component mixture for $q
   = 1, 2, \ldots$}$ when mixtures are formed from any smooth
 families, including discrete, continuous and multivariate densities;
 and (3) testing for the presence of a signal when the probability
 density functions under the null and alternative hypotheses belong to
 different parametric families which occurs in physics applications
 \citep{pilla:05}.

\subsection{Inference in Spatial Scan Statistics}  

In the scan statistics problem, one observes a random field (such as a
point process) in a region of interest. The goal is to detect unusual
behavior in subregions, where the behavior of the field differs
significantly from the background. Applications include mammography;
automatic target recognition; disease clustering and minefield
detection.

In the classical formulation of the scan statistic (see \cite{glaz:01}
and the references therein), a rectangular window is scanned across
the data, with high values of the statistic indicating a local
departure from uniformity. In contrast, the methods developed in this
article are applicable to smooth scanning processes, where the window
is tapered, rather than having sharp boundaries. The null density
\(f(\cdot; \mbl)\) represents the background model while the scan
window \(\psi(\cdot; \mbt)\) represents departure from the background
at location \(\mbt\).

\subsection{Main Results}

We create a general family of models referred to as perturbation
models that encompass a large class of statistical problems. Our
treatment of the nuisance parameters under the null hypothesis is
quite general. The inferential theory developed in this article
provides a solution to an important class of statistical problems
involving loss of identifiability and/or when some of the parameters
are on the boundary of the parametric space.  The main contributions
of this article are as follows.
\begin{enumerate}
\vspace{-0.15in}
\item In Section \ref{sec-score}, we propose a novel test statistic
  based on the score process, denoted by $\mT$, for detecting the
  presence of perturbation and derive its fundamental properties. In
  particular, it is shown that the test statistic $\mT$ based on the
  score process is asymptotically equivalent to the LRT statistic.
  \vspace{-0.1in}
\item In Section \ref{sec-perturb}, we derive a general inferential
  theory for approximating the asymptotic null distribution of $\mT$. It
  is shown that the asymptotic distribution of $\mT$ under $\H_0$ equals
  $\sup_{\mbt} \, Z(\mbt)$, where $Z(\mbt)$ is a differentiable
  Gaussian random field with continuous sample paths. Therefore, the
  goal becomes finding approximations for $\mP(\sup_{\mbt} \, Z(\mbt)
  \geq c)$ for any large $c \in \R$ in order to determine the
  quantiles of $\mT$. As eloquently pointed out by \cite{adler:00}, this
  problem occurs in a large number of different applications including
  in image processing \citep{worsley1:95}. We describe a connection
  between $Z(\mbt)$ and a differentiable manifold (curve, surface,
  etc.)  through the {\em Karhunen-Lo\`{e}ve expansion}. The
  Karhunen-L\`{o}eve expansion converts the high-dimensional Gaussian
  probability problem into that of a chi-squared random variable and
  uniformly distributed random variables over the surfaces of spheres
  \citep{adler:00}.  \vspace{-0.1in}
\item Our technique is based on the long-established and elegant
  geometric result known as the {\em volume-of-tube formula}
  \citep{hotelling:39,weyl:39,naiman:90}. The problem of evaluating
  the Gaussian random field significance probabilities (i.e., tail
  probability for the asymptotic null distribution of $\mT$) for testing
  the hypothesis (\ref{eq:bashyp}) is reduced to that of determining
  the volume-of-tube about a manifold on the surface of a hypersphere
  (see Section \ref{sec-distr}). The novelty here lies in deriving
  explicit expressions for the geometric constants appearing in the
  volume-of-tube formula with boundaries; consequently, one can
  approximate the quantiles of the statistic $\mT$ for detecting the
  presence of perturbation. We also address the difficult and yet
  important problem of presence of singularities in the score process.
  \vspace{-0.1in}
\item In Section \ref{sec-nuisance}, the results of Section
  \ref{sec-perturb} are extended to the case where the null density is
  characterized by a vector of nuisance parameters.  \vspace{-0.1in}
\item An age old and fundamental question of determining the order of
  a mixture model is solved in Section \ref{sec-mixture}. In
  particular, building on the perturbation theory, we develop
  inferential methods for approximating the quantiles of the test
  statistic for determining the mixture complexity. The flexibility
  and general applicability of the methodology is demonstrated through
  univariate and multivariate mixture families. Furthermore, it is
  shown that the results of \cite{lindsay:95}, \cite{lin2:97} and
  \citet{chen1:01} become special cases of our general and broadly
  applicable theory.
\end{enumerate}

The paper concludes with a discussion of the relative merits of the
perturbation theory in Section \ref{sec-discuss}. In Section
\ref{sec-proof}, we derive the proofs of our general results. Explicit
expressions for the geometric constants that appear in the volume-of-tube 
formula are derived in Appendix A. 

\setcounter{equation}{0} \def\theequation{\thesection.\arabic{equation}}

\section{A Score Process and its Fundamental Properties} 
\label{sec-score}
 
In this section, we derive a score process and its fundamental
properties that are required for the testing problem
(\ref{eq:bashyp}). As a first step, we assume that $\mbl$ is fixed or
known so that \(f(; \mbl)\) is completely specified and the density
(\ref{eq:pbm}) can be expressed simply as $p(; \eta, \mbt)$; 
however, theory for the general case of an unknown $\mbl$ will be derived 
in Section \ref{sec-nuisance}. 

\subsection{Loglikelihood Ratio Process} 
If $\mbt$ is fixed at a particular value, then the testing problem
(\ref{eq:bashyp}) becomes routine. However, the nuisance parameter
vector $\mbt$ can assume any value under $\H_0$; therefore, the testing
problem is non-regular. The loglikelihood function based on the
perturbation model (\ref{eq:pbm}) is $l(\eta, \mbt|\bx) = \sum_{i =
  1}^n \log \, [(1 - \eta) \, f(x_i; \mbl) + \eta \, \psi(x_i;
\mbt)]$. For a fixed \(\mbt\), $l(\eta, \mbt|\bx)$ is a
concave function of \(\eta\) and hence there exists a unique maximizer
\(\we_{\mbt} \in [0, 1]\). In general, there is no closed form
solution for \(\we_{\mbt}\); however, the estimator can be found as a
solution to
\begin{equation}
  \label{eq:eps}
  \sum_{i = 1}^n \frac{ [\psi(x_i; \mbt) - f(x_i; \mbl)]
  }{p(x_i; \eta, \mbt)} = 0 
\end{equation}
if a solution in \((0, 1)\) exists; otherwise the estimator will be at
one of the end-points.  This leads to a corresponding {\em
loglikelihood ratio process} $l^{\star}(\mbt|\bx) = l(\we_{\mbt},
\mbt|\bx) - l(0, 0|\bx)$. Considered as a function of \(\mbt\), the
process \(l^{\star}(\mbt|\bx)\) may be used as a diagnostic tool, with
large values indicating the presence of perturbation. The maximum
likelihood estimator (MLE) of \(\mbt\) is the maximizer of
\(l^{\star}(\mbt|\bx)\). However, maximizing this process is 
computationally intensive, since \(l^{\star}(\mbt|\bx)\) may have many
local maxima. Any strategy for finding the global maximum has to
involve an exhaustive search, which in turn requires solving
(\ref{eq:eps}) for each fixed \(\mbt\). In the next section, we
derive an alternative technique that will combat these difficulties.

\subsection{The Score Process: Theory}

In this section, we propose a novel technique based on a {\em
score process} defined as
\begin{eqnarray}
  \label{eq:scorep}
  S(\mbt) &:=& \left. \frac{\partial}{\partial \eta}
  l(\eta, \mbt, \mbl|\bx) \right|_{\eta = 0} = \sum_{i = 1}^n \left[
  \frac{\psi(x_i; \mbt)}{f(x_i; \mbl)} - 1 \right].
\end{eqnarray}
The interest is in the parameter vector $\mbt$ and since $\mbl$ is
fixed for now, for exposition, we drop $\mbl$ from the expressions and
simply write $S(\mbt), S^{\star}(\mbt), Z(\mbt)$, etc. 

The score process has several elegant features: (1) it is not as
computationally intensive as the likelihood ratio process and (2) its
explicit representation makes statistical inference tractable.  It is
shown in Theorem \ref{th:scorethm} (below) that the score process has
mean zero when there is no perturbation (i.e., \(\eta = 0\)) and
$\mE[S(\mbt)] > 0$ when there is a perturbation at \(\mbt = \mbt_0\),
the true parameter vector. This suggests that peaks in the score
process provide evidence for the presence of perturbation. However,
$S(\mbt)$ can exhibit high random variability and the variance may
have substantial dependence on \(\mbt\). To combat this difficulty, we
propose the {\em normalized score process} defined as 
\begin{eqnarray*}
  S^{\star}(\mbt) := \frac{S(\mbt)}{\sqrt{n \, \mC(\mbt, \mbt)}}, 
\end{eqnarray*}
where the {\em covariance function} is defined as 
\begin{eqnarray}
\label{eq:covfn}
  \mC(\mbt, \mbt^{\dag}) &:=& \int \frac{ [\psi(x; \mbt) - f(x; \mbl)]
  \, [\psi(x; \mbt^{\dag}) - f(x; \mbl)]}{f(x; \mbl)} \, dx \nonumber \\
  &=& \int \frac{\psi(x; \mbt) \, \psi(x; \mbt^{\dag})}{f(x; \mbl)}
  \, dx - 1.
\end{eqnarray}
The covariance function $\mC(\mbt, \mbt^{\dag})$ has an analytical
expression for certain choices of $f(\cdot; \mbl)$ and $\psi(\cdot;
\mbt)$ while in other cases numerical integration is required. 

The following conditions are assumed for deriving the large-sample
theory.

\assumption 
\label{as:pspace}
The parameter space $\mbT$ is a compact and a convex subset of $\R^d$
for some integer $d$. 

\assumption 
\label{as:finitev}
The covariance function satisfies $\mC(\mbt, \mbt) < \infty$ for all
$\mbt \in \mbT$. 

\assumption 
\label{as:supp}
For each \(\mbt \in \mbT\), \( \supp[\psi(\,\cdot\,; \mbt)]
\subset \supp[f(\,\cdot\,; \mbl)] \), where `$\supp$' refers to
the support of a density.

In the following theorem, we characterize some fundamental properties
of the score and normalized score processes.
\begin{theorem} 
\label{th:scorethm}
Suppose assumptions \ref{as:finitev} and \ref{as:supp} hold: (1) Under $\H_0$, the score 
process has mean $\mE[S(\mbt)] = 0 \; \mbox{for all} \; \mbt$ with a
covariance function $\, {\rm cov}[S(\mbt), S(\mbt^{\dag})] = n \,
\mC(\mbt, \mbt^{\dag})$, where $\mC(\mbt, \mbt^{\dag})$ is defined in
(\ref{eq:covfn}); (2) under $\H_1$,
\begin{equation}
  \label{eq:altmean}
  \mE[S(\mbt)] = {\it n} \, \eta \, \mC(\mbt, \mbt_0);
\end{equation}
and (3) under $\H_1$, the expectation of the normalized score
process is 
\begin{equation}
\label{eq:ansmean} 
  \mE[S^{\star}(\mbt)] = \sqrt{\it n} \, \eta \, \frac{ \mC(\mbt,
  \mbt_0) }{\sqrt{\mC(\mbt, \mbt)}} \le \eta \sqrt{\it n 
  \, \mC(\mbt_0, \mbt_0)}
\end{equation}
with equality at \(\mbt = \mbt_0\).
\end{theorem}
{\em Proof.} Under $\H_1$, it follows that 
\begin{eqnarray}
  \mE[S(\mbt)] &=& n \int \left[ \frac{ \psi(x; \mbt)}{f(x;
  \mbl)} - 1\right] p(x; \eta, \mbt_0) \, dx \nonumber \\
  &=& n \, \eta \int   \frac{ \psi(x; \mbt) \,
  \psi(x; \mbt^{\dag})}{f(x; \mbl)} \, dx = n \, \eta \, 
  \mC(\mbt, \mbt_0) \nonumber
\end{eqnarray}
which yields the result (\ref{eq:altmean}). Similarly, one can derive
the mean and covariance functions in part 1 of the theorem. The bound
(\ref{eq:ansmean}) is established by noting that \(\mC(\mbt,
\mbt_0)\) is a covariance function and therefore satisfies the
Cauchy-Schwartz inequality $\mC(\mbt, \mbt_0) \le \sqrt{ \mC(\mbt,
\mbt) \, \mC(\mbt_0, \mbt_0) }$. \hfill \rule{2mm}{2mm}

The motivation for using the score processes lies in part 3 of Theorem
\ref{th:scorethm}: The expectation of $S^{\star}(\mbt)$ is
maximized at $\mbt_0$. Therefore, the supremum of the process
$S^{\star}(\mbt)$ can serve as a test statistic for the hypothesis
(\ref{eq:bashyp}). If $\H_0$ is rejected, then the maximizer of
$S^{\star}(\mbt)$ serves as a point estimator of \(\mbt\).  The
final result of this section establishes the asymptotic equivalence
between the score and loglikelihood processes; the proof is given in
Section \ref{sec-proof}.

\begin{theorem} 
\label{th:lrtscore}
The score process and loglikelihood ratio process are asymptotically
equivalent, in the sense that $l^{\star}(\mbt|\bx) =
\frac{1}{2}[\max\{0, S^{\star}(\mbt) \}]^2 + o_p(1)$ as $n \to
\infty$.   
\end{theorem}

\setcounter{equation}{0}
\def\theequation{\thesection.\arabic{equation}}

\section{Testing for the Presence of Perturbation}
\label{sec-perturb}

We first propose a statistic for the testing problem (\ref{eq:bashyp})
and next derive its asymptotic null distribution. From the motivation
presented in the previous section, it is natural to define a
statistic for testing the hypothesis (\ref{eq:bashyp}) as 
\begin{eqnarray}
  \label{eq:ts}
  \mT := \supt \, S^{\star}(\mbt).  
\end{eqnarray} 
Except in special cases, the distribution of $\mT$ cannot be expressed
analytically. Our next goal is to derive an asymptotic distribution of
\(\mT\) under \(\H_0\) for determining approximate quantiles of the test
statistic. As a first step, we establish that under $\H_0$ the
distribution of \(\mT\) is asymptotically equivalent to the distribution
of the supremum of a Gaussian random field. Next, we derive
approximations for the tail probability of the supremum of a Gaussian
random field using the Karhunen-Lo\`{e}ve expansion and the
volume-of-tube formula.

The volume-of-tube problem for curves (i.e., $d = 1$) was first
studied by \cite{hotelling:39} in the context of significance testing
for nonlinear regression. In a second pioneering paper, \cite{weyl:39}
extended the work of Hotelling to higher-dimensional manifolds (i.e.,
$d \geq 2$), deriving elegant expressions for the volume-of-tube of
manifolds lying in a hypersphere. \cite{naiman:90} further extended
the Hotelling-Weyl results to cases where the manifold has boundaries.
\cite{suntube} studied higher order terms for Gaussian processes and
fields. Important statistical problems to which the volume-of-tube
formula has been applied include non-linear regression
\citep{hotelling:39,knowles:89}, projection pursuit 
\citep{johansen:90}, testing for multinomial mixture models
\citep{lindsay:95,lin2:97}, simultaneous confidence bands
[\citet{naiman:87}, \citet{sunloascb} and Chapter 9 of
\citet{loader:99}] and inference under convex cone alternatives for
correlated data \citep{pilla:06}.

The following assumptions are required for the development of
inferential theory.    

\assumption 
\label{as:diffpsi}
  For all $x \in \X$, the perturbation density \(\psi(x; \mbt)\) is a 
  twice differentiable, while
\begin{eqnarray*}
  \int \frac{\psi^{\prime}(x, \mbt)^2}{f(x, \mbl)} dx < \infty \quad
  \mbox{and} \quad \int \frac{\psi^{\prime \prime}(x, \mbt)^2}{f(x,
  \mbl)} dx < \infty,
\end{eqnarray*}
where $\prime$ denotes differentiation with respect to
\(\mbt \in \mbT\). In the multi-parameter case, all first and
second-order partial derivatives are assumed to satisfy the
integrability condition as well.

\assumption
\label{as:poscov}
The covariance function $\mC(\mbt, \mbt)$ is positive in $\mbt$;
equivalently, \(f(\cdot; \mbl)\) is not identically equal to
\(\psi(\cdot; \mbt)\) for any \(\mbt \in \mbT\).    

The assumption \ref{as:poscov} fails in several important problems
including mixture models, leading to singularities in the score
process. In Section \ref{sec-singular}, we derive modifications to our
theory to handle this difficult but important problem.

Let $\{Z(\mbt)\!: \mbt \in \mbT \subset \R^d\}$ be a $d$-dimensional
differentiable Gaussian random field with continuous sample paths,
with mean zero and covariance function 
\begin{eqnarray}
  \label{eq:rho}
  \rho(\mbt, \mbt^{\dag}) := \mE\left[ Z(\mbt) Z(\mbt^{\dag}) \right]
  = \frac{ \mC(\mbt, \mbt^{\dag}) }{ \sqrt{ \mC(\mbt, \mbt) \, 
    \mC(\mbt^{\dag}, \mbt^{\dag}) } }.
\end{eqnarray}
Under assumptions \ref{as:diffpsi} and \ref{as:poscov}, the
asymptotic null distribution of \(\mT\) is the supremum of a 
Gaussian random field, expressed explicitly as
\begin{eqnarray*} 
  Z(\mbt) = [\mC(\mbt, \mbt)]^{-1/2}  \int \left[ \frac{
  \psi(x; \mbt)}{f(x; \mbl)} - 1 \right] \sqrt{f(x, \mbl)} \; W(d x), 
\end{eqnarray*}
where \(W\) is the standard {\em Brownian sheet}. 

\begin{theorem} \label{th:maxsstar}
Suppose that assumptions \ref{as:pspace} to \ref{as:poscov}
hold. Under  \(\H_0\),
\begin{equation}
  \label{eq:ascrv} 
  \mP \left(\mT \ge c \right) \longrightarrow \mP \left(\supt
  Z(\mbt) \ge c \right) \quad \mbox{as} \quad n \rightarrow 
  \infty \quad \mbox{for any} \quad c \in \R.
\end{equation}
\end{theorem}

Theorem \ref{th:maxsstar} will be proved in Section
\ref{sec-proof}. Generally, there is no exact result for finding
$\mP(\sup_{\mbt} Z(\mbt) \ge c)$ \citep{adler:00}. The result of
Theorem \ref{th:maxsstar} holds even if we relax assumption
\ref{as:diffpsi}. Our proof relies only on the assumption of first
derivative of $\psi(\cdot; \mbt)$; however, the second derivative
conditions are required for the explicit probability approximations
derived later using the volume-of-tube-formula.

The problem of approximating the distribution of the supremum of a
smooth Gaussian random field (i.e., finding $\mP(\sup_{\mbt} Z(\mbt)
\ge c)$  for large $c$) can be addressed using several different
techniques: (1) methods based on the Hotelling-Weyl
\citep{hotelling:39,weyl:39} volume-of-tube formula with boundary
corrections \citep{naiman:90}; (2) expected Euler characteristic
methods \citep{siegmund:95,worsley:01}; (3) approaches based on
counting the local maxima and upcrossings; and (4) Rice formula
\citep{siegmund:93,azais:05}.  All these techniques lead to similar
results for practical purposes (see \cite{adler:00} for
discussion). Some formal equivalence results between the tube formula
and the expected Euler characteristic methods have been derived by
\cite{take:02}. In this article, for the development of 
inferential theory for perturbation models, we adopt the
volume-of-tube formula technique for its relatively simple geometric
interpretation and the flexibility to yield explicit results for
higher-order boundary corrections. The disadvantage of the tube
approach is that it is directly applicable only to processes that are
Gaussian or Gaussian-like \citep{adler:00}. 

\subsection{The Karhunen-Lo\`{e}ve Expansion} 

In this section, we construct a sequence of finite-dimensional
approximation to the Gaussian random field \(Z(\mbt)\) using the
Karhunen-Lo\`{e}ve expansion. Although Karhunen-Lo\`{e}ve expansion is
most convenient, any other uniformly convergent approximation, such as
a cubic spline interpolant on a grid of \(\mbT\) is also applicable.

While some of the core ideas in this section are known, there does not
exist a complete statement of the results in the form that are
required for the general testing problem (\ref{eq:bashyp}). In
particular, addressing the following scenarios are of fundamental
importance: (1) \(\mbT\) is a hyper-rectangle or a similar polygonal
region with boundaries of various orders (edges, corners and so on)
and (2) the score process $S(\mbt)$ has singularities.

A concise presentation of the Karhunen-Lo\`{e}ve expansion can be
found in Section III.3 of \citet{adler:90}. The Karhunen-Lo\`{e}ve
expansion of $Z(\mbt)$ is the \emph{uniformly convergent} series
expansion
\begin{equation}
  \label{eq:klexp}
  Z(\mbt) = \sum_{k = 1}^{\infty} \mz_k \; \xi_k(\mbt) 
    = \ip{\mbu}{\mbx(\mbt)},
\end{equation}
where $\mz_k$ is an i.i.d.\ standard Gaussian random variable,
\(\{\xi_k(\mbt)\}_{k = 1}^{\infty}\) is a sequence of twice
continuously differentiable functions, while $\mbu$ and \(\mbx(\mbt)\)
are the corresponding vector counterparts. The covariance function
(\ref{eq:rho}) can be explicitly expressed as \begin{eqnarray}
  \label{eq:rho1}
  \rho(\mbt, \mbt^{\dag}) &=& \sum_{k = 1}^{\infty}
    \xi_k(\mbt) \, \xi_k(\mbt^{\dag})
\end{eqnarray}
and $\mz_k = \mu_k^{-1} \int_{\mbT} \xi_k(\mbt) \, Z(\mbt)
\, d \mbt$, where \(\mu_k = \int_{\mbt} \xi_k^2 (\mbt) \, d \mbt\).

It is necessary for $Z(\mbt)$ to have a finite Karhunen-Lo\`{e}ve
expansion for the application of the volume-of-tube formula. When the
expansion is infinite, the series is truncated at \(J\) terms to yield
\begin{eqnarray}
  Z_J(\mbt) &:=& \sum_{k = 1}^{J - 1} \mz_k \, \xi_k(\mbt)
    + \mz_0 \, \sqrt{\sum_{k = J + 1}^{\infty} [\xi_k(\mbt)]^2}
    = \ip{\mbu_J}{\mbx_J(\mbt)},
 \label{eq:truncz}
\end{eqnarray}
where \(\mz_0 \sim N(0, 1)\) and is independent of \(\mz_1, \mz_2,
\ldots \), \(\mbu_J = \left (\mz_0, \ldots, \mz_{J - 1} \right
)^T\) and \(\mbx_J(\mbt)\) is the corresponding truncated version of
the sequence \(\{\xi_k(\mbt)\}_{k = 1}^\infty\). The covariance
function of \(Z_J(\mbt)\) can be expressed as 
\begin{eqnarray}
  \label{eq:rhoJ}
  \rho_J(\mbt, \mbt^{\dag}) = \sum_{k = 1}^{J - 1} \xi_k(\mbt) \;
  \xi_k(\mbt^{\dag}) + \sqrt{ \sum_{k = J}^{\infty} \xi_k^2(\mbt)
  \sum_{k = J}^{\infty} \xi_k^2(\mbt^{\dag})} =
  \ip{\mbx_J(\mbt)}{\mbx_J(\mbt^{\dag}) }. 
\end{eqnarray}
The final term in (\ref{eq:truncz}) has been chosen to preserve
unit variance; i.e., $\mV[Z_J(\mbt)] = \rho_J(\mbt, \mbt) = 1$.  

\subsection{Distribution of the Supremum of $Z(\mbt)$}
\label{sec-distr}

In this section, we provide an approximation to $\sup_{\mbt} Z(\mbt)$
under a very general assumption that $\M$ is a manifold with a
piecewise smooth boundary. This result, combined with Theorem
\ref{th:maxsstar} provides an elegant approximation to the asymptotic
null distribution of the test statistic \(\mT\). The primary goal is
to approximate the asymptotic probability in (\ref{eq:ascrv}) when $c
\in \R$ is large, $\mbt \in \mbT \subset \R^d$ and $d \ge 1$. 

Conditioning on the length of the vector $\mbu_J$, 
\begin{eqnarray}
  \mP \left(\supt Z_J(\mbt) \ge c \right)
    &=& \mP\left(\supt \, \ip{\mbu_J}{\mbx_J(\mbt)} \ge c \right)
    \nonumber \\
   &=& \mP \left(\supt \, \ip{\frac{\mbu_J}{\|\mbu_J\|}}{\mbx_J(\mbt)}
    \ge \frac{c}{\|\mbu_J\|} \right) \nonumber \\
    &=& \int_{c^2}^{\infty} \mP\left(\supt \, 
    \ip{\bU_J}{\mbx_J(\mbt)} \ge \frac{c}{\sqrt{y}} \right)
    h_J(y) \, dy, \label{eq:iprob}
\end{eqnarray}
where the $J$-dimensional random vector $\bU_J = (\mz_0/\|\mbu_J\|,
\ldots, \mz_{J - 1}/\|\mbu_J\|)^T$ is uniformly distributed on the
unit sphere $\S^{(J - 1)}$ embedded in $\R^J$, $\mbx(\mbt)$ is a curve
in $\S^{(J - 1)}$ and \(h_J(y)\) is the $\chi^2$ density with \(J\)
degrees of freedom. Consequently, the goal becomes evaluating the
distribution of the supremum of a uniform process in
(\ref{eq:iprob}). 

First, note that the inner product \(\ip{\bU_J}{\mbx_J(\mbt)}\) is
bounded by 1 (using the Cauchy-Schwarz inequality) enabling the
restriction of \(c/\sqrt{y} < 1\) or \(c^2 < y < \infty\). Since $\|
\bU_J - \mbx_J(\mbt) \|^2 = \|\bU_J \|^2 + \| \mbx_J(\mbt) \|^2 - 2
\ip{\bU_J}{\mbx_J(\mbt)} = 2[1 - \ip{\bU_J}{\mbx_J(\mbt)}]$, it
follows that, for any \(w \in (0, 1)\), \(\ip{\bU_J}{\mbx_J(\mbt)} \ge
w\) if and only if \( \| \bU_J - \mbx_J(\mbt) \| \le r := \sqrt{2 (1 -
  w)}\). Therefore,
\begin{eqnarray}
  \mP \left(\supt \, \ip{\bU_J}{\mbx_J(\mbt)} \ge w \right)
  &=& \mP \left(\inft \|\bU_J - \mbx_J(\mbt)\| \le r \right)
  \nonumber\\  
  &=& \mP[\bU_J \in \T(r, \M)] = \frac{ \vartheta(r, \M)}{A_J},
  \label{eq:tuberat}
\end{eqnarray}
where $\vartheta(r, \M)$ denotes the volume of \(\T(r, \M)\)---a
\emph{tube} of radius $r$ around the \emph{manifold} $\M := \{
\mbx_J(\mbt)\!: \mbt \in \mbT \subset \R^d\}$, and \(A_J = 2
\pi^{J/2}/\Gamma(J/2)\) is the $(J - 1)$-dimensional volume of the
unit sphere $\S^{(J - 1)}$. The last expression follows since
\(\bU_J\) is uniformly distributed over $\S^{(J - 1)}$. 

\begin{figure}[htb]
\centerline{\scalebox{0.6}{\includegraphics{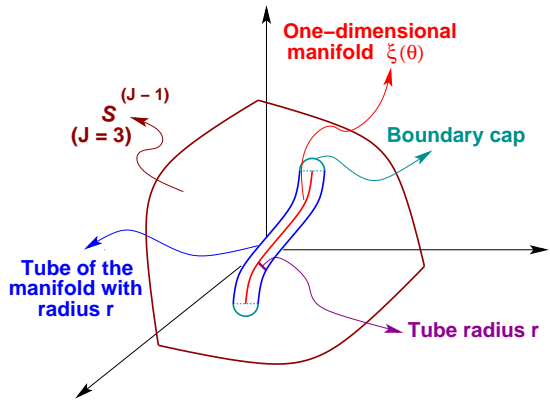}}}
\caption{Tube of radius $r$ around a one-dimensional manifold (curve) 
with boundaries embedded in \(\S^2\).}
\label{fig:tube}
\end{figure}

\remark Finding the distribution of the supremum of a Gaussian random
field $Z(\mbt)$ is now reduced to that of determining the
volume-of-tube of the manifold \(\M\). The solution to this problem
depends on the geometry of \(\M\). When the set \(\mbT\) is
one-dimensional (i.e., $d = 1$) and \(\mbx_J(\mbt)\) is continuous,
then $\M$ is a curve on the unit sphere $\S^1$ and the tube consists
of a main ``cylindrical'' section plus the two boundary caps as shown
in Fig.\ \ref{fig:tube}. In this case, results of \cite{hotelling:39}
and \cite{naiman:90} yield the approximation \begin{eqnarray*}
\vartheta(r, \M) \approx \k_0 \, \frac{A_J}{A_2} \mP \left[B_{1, (J - 2)/2} 
\ge w^2 \right] + \, \ell_0 \, \frac{A_J}{2 A_1} \mP \left[ B_{1/2, (J - 1)/2} 
  \ge w^2 \right], \end{eqnarray*} where \(\k_0\) is the length of the
  manifold \(\M\), $B_{a, b}$ is the beta density with parameters $a$
  and $b$ and \(\ell_0 = 2\) is the number of end-points. Introducing
  \(\ell_0\) allows us to treat cases where \(\M\) consists of two or
  more disconnected segments (due to singularities in the score
  process), which is a common phenomena in the context of mixture
  models. The volume-of-tube formula is exact whenever $r$ is less
  than a critical radius \(r_0\) (equivalently, $w_0 \leq w \leq 1$)
  which depends on the curvature of \(\M\).

Application of the volume-of-tube formula to a Gaussian random field
leads to the main result of this section.

\begin{theorem}
\label{th:sig3d}
Under assumptions \ref{as:pspace} to \ref{as:poscov}, the
distribution of $\sup_{\mbt} Z(\mbt)$ for a general $d$ is given by
\begin{eqnarray}
   \mP \left(\supt Z(\mbt) \ge c \right) = \sum_{t = 0}^d
   \frac{\zeta_t}{A_{d + 1 - t}} \mP\left( \chi^2_{d + 1 - t} \ge
   c^2\right) + o[c^{-1} \exp(-c^2/2)], \label{eq:tubechis}
\end{eqnarray}
as $c \to \infty$, where \(A_t = 2 \pi^{t/2}/\Gamma(t/2)\) is the $(t
- 1)$-dimensional volume of $\S^{(t - 1)}$ in \(\R^t\) and \(\zeta_t\)
are the geometric constants derived in Appendix A.
\end{theorem}

\noindent {\em Multinomial Mixture Problem:} Equation (4.19) of
\cite{lindsay:95}, derived in the context of multinomial mixture
models, is a special case of Theorem \ref{th:sig3d} (see also
\cite{lin2:97} for bounds). This connection is explored further in
Section \ref{sec-mixture}. It is important to note that for
multinomial mixture models, the Karhunen-Lo{\`e}ve expansion is
finite.

\remark Although the proof of Theorem \ref{th:sig3d}, derived in
Section \ref{sec-proof}, uses the Karhunen-Lo\`{e}ve expansion, it is
not necessary to find this expansion since one can determine the
geometric constants $\zeta_t$s appearing in (\ref{eq:tubechis})
entirely from the covariance function \(\mC(\mbt,
\mbt^{\dag})\). However, it is necessary to consider the geometry
of the manifold \(\M\) in order to treat the boundary corrections,
particularly when $d > 2$.

\subsection{Singularities in the Score Process}
\label{sec-singular}

One of the conditions required for Theorem \ref{th:sig3d} is that
$\mC(\mbt, \mbt)$ is positive for all \(\mbt \in \mbT\). This
condition is violated when \(f(\cdot; \mbl) = \psi(\cdot; \mbt)\) for
some \(\mbt\). This is a commonly occurring phenomena in the context
of finite mixture models. Therefore, we need to consider more
carefully the behavior of the score process near \(\mbt =
\mbt_0\). Let $S^{\prime}(\mbt) = \partial S(\mbt)/\partial \mbt$ and
$\mV[S^{\prime}(\mbt)]$ be the variance of $S^{\prime}(\mbt)$ so that
$S(\mbt) = (\mbt - \mbt_0) S^{\prime}(\mbt_0) + o(\mbt - \mbt_0)$, $n
\, \mC(\mbt, \mbt) = (\mbt - \mbt_0)^2 \, \mV[S^{\prime}(\mbt_0)] +
o[(\mbt - \mbt_0)^2]$ and $S^{\star}(\mbt) = \sgn(\mbt - \mbt_0)
S^{\prime}(\mbt_0)/\sqrt{\mV[S^{\prime}(\mbt_0)]} + o(\mbt
- \mbt_0)$, where `$\sgn$' is the sign function. In particular, this
implies that the process ``flips'' and 
\begin{equation}
  \lim_{\mbt \to \mbt_0^-} S^{\star}(\mbt) =
    - \lim_{\mbt \to \mbt_0^+} S^{\star}(\mbt).
    \label{eq:flip} 
\end{equation}
Correspondingly, \(\mbx(\mbt_0^-) = -\mbx(\mbt_0^+)\). In effect,
the manifold \(\M\) has two pieces and four boundary points. The
result in Theorem \ref{th:sig3d} still holds; however, $\ell_0 = 4$. 

\setcounter{equation}{0}
\def\theequation{\thesection.\arabic{equation}}

\section{Nuisance Parameters under the Null Model}
\label{sec-nuisance}

In this section, we derive general theory for the case of unknown
nuisance parameter vector $\mbl$. We derive a series of fundamental
results that provide a ``linearization of the score process'' (defined
below) to identify the correct covariance function (see Theorem
\ref{th:nucov} below) for this setting. We replace $\mbl$ by
$\widehat{\mbl}$, the MLE of $\mbl$, and assume that the MLE satisfies
the necessary regularity conditions stated by \cite{chernoff:54}. Our
goal is to find an appropriate normalizing factor for the score
process and in turn apply the volume-of-tube formula for approximating
the asymptotic null distribution of $\mT$.

In the context of finite mixture models, the null density \(f(\cdot;
\mbl)\) is equivalent to the mixture density ${\rm g}(\cdot; \Q_m)$
representing an \(m\)-component mixture model with \(\mbl \equiv
\Q_m\) containing a vector of support points and the corresponding
mixing weights. The score process is searching for an \((m +
1)\)st component.

If \(\mbl\) is estimated via the ML method, then under \(\H_0\), the
score process can be expressed as 
\begin{eqnarray*}
  S(\mbt|\widehat{\mbl}) := \sum_{i = 1}^n \left[\frac{\psi(x_i;
  \mbt)}{f(x_i; \widehat{\mbl})} - 1 \right]. 
\end{eqnarray*}
The statistic $\mT$ will still be the supremum (over \(\mbt\))
of the normalized score process; however, estimating the nuisance
parameter vector \(\mbl\) means that the covariance function
\(\mC(\mbt, \mbt^{\dag})\) defined in (\ref{eq:covfn}) is no longer
appropriate for normalizing the score process. 

As a first step, it is assumed that the MLE \(\widehat{\mbl}\) under
$\H_0$ satisfies the required conditions for the second-order
asymptotic theory \citep{lehmann:99}. Hence, the following results
hold: 
\begin{eqnarray} 
  (\widehat{\mbl} - \mbl_0) &=& \left[n {\bf I}(\mbl_0) \right]^{-1}
  \sum_{i = 1}^n \nabla \, l(\mbl|x_i) + o_p(n^{-1/2}) 
     \label{eq:asthhat} 
\end{eqnarray} 
and $n^{-1/2} \sum_{i = 1}^n \nabla \, l(\mbl_0|x_i) \toinlaw
N[0, {\bf I}(\mbl_0)]$ as $n \to \infty$, where \(\mbl_0\) is the
true null parameter vector, $\toinlaw$ indicates convergence in
distribution, \({\bf I}(\mbl_0)\) is the Fisher information
matrix and \(\nabla \, l(\mbl|x)\) is the vector of partial
derivatives of \(l(\mbl|x) = \log f(x; \mbl)\) with respect to
\(\mbl\).

\begin{theorem} Suppose that assumptions \ref{as:finitev} to
  \ref{as:poscov} hold. Under $\H_0$ with the true null parameter
  vector \(\mbl_0\), the score process has the asymptotic representation of
\begin{eqnarray*}
   S(\mbt| \widehat{\mbl}) = S(\mbt| \mbl_0) - \mC^T(\mbt|\mbl_0) \,
   [{\bf I}(\mbl_0)]^{-1} \, \sum_{i = 1}^n \nabla \, l(\mbl_0|x_i) + 
   o_p(n^{1/2}),  
\end{eqnarray*}
where \(o_p(n^{1/2})\) is uniform in \(\mbt\) and \(\mC(\mbt|\mbl_0)\)
is the covariance vector defined as  
\begin{eqnarray*}
  \mC(\mbt| \mbl_0) := \cov \left[ \left(\frac{
    \psi(x_1; \mbt)}{f(x_1; \mbl_0)} - 1 \right),
    \, \nabla \, l(\mbl_0|x_1) \right] = \int \psi(x;
    \mbt) \, \nabla \, l(\mbl_0|x) \, dx. 
\end{eqnarray*}
\end{theorem}
{\em Proof.} By expanding the score process in a Taylor series around
\(\mbl_0\), we obtain
\begin{eqnarray*}
  S(\mbt|\widehat{\mbl}) = S(\mbt|\mbl_0)
    + (\widehat{\mbl} - \mbl_0)^T 
       \left . \frac{ \partial }{\partial \mbl} S(\mbt|\mbl)
      \right|_{\mbl = \tilde\mbl},
\end{eqnarray*}
where \(\widetilde\mbl \in [\mbl_0, \widehat{\mbl}]\). Direct
calculation shows that 
\[
  \left. \frac{1}{n} \frac{\partial}{\partial \mbl} S(\mbt|\mbl) 
      \right|_{\mbl = \widetilde\mbl}
  = \; - \frac{1}{n}  \sum_{i = 1}^n \frac{ \psi(x_i; \mbt)
    }{f(x_i; \mbl_0)} \, \nabla \, l(\mbl|x_i).
\]
From the uniform strong law of large numbers and the fact that
\(\widetilde\lambda \toinas \lambda_0\), it follows that  
\begin{eqnarray*}
  \left . \frac{1}{n} \frac{\partial}{\partial \mbl} S(\mbt|\mbl)
    \right|_{\mbl = \widetilde\mbl}
    &\toinas& - \, \cov \left[ \frac{\psi(x_1; \mbt)
      }{f(x_1; \mbl_0)} - 1, \, \nabla \, l(\mbl_0|x_1) \right] \\ 
     &=& - \, \mC(\mbt|\mbl_0) \quad \mbox{as} \quad n \to \infty.
\end{eqnarray*}
It follows from assumption \ref{as:pspace} and the continuity of
$\mbT$ that the convergence is uniform in \(\mbt\). Combining this
result with (\ref{eq:asthhat}) completes the proof.  \hfill
\rule{2mm}{2mm} 

\begin{theorem}\label{th:nucov}
The process
\begin{equation}
  \label{eq:slth0}
    n^{-1/2} \, S(\mbt| \mbl_0) - \, n^{-1/2} \, \mC^T(\mbt| \mbl_0)
    \; [{\bf I}(\mbl_0)]^{-1} \, \sum_{i = 1}^n \nabla \, l(\mbl_0|x_i)
\end{equation}
has the covariance function
\begin{eqnarray}
  \label{eq:covfn2} 
  \mC^{\star}(\mbt, \mbt^{\dag}) =
    \mC(\mbt, \mbt^{\dag}) - [\mC^T(\mbt|\mbl_0)] \; 
    [{\bf I}(\mbl_0)]^{-1} \, \mC(\mbt^{\dag}|\mbl_0),
\end{eqnarray}
where \(\mC(\mbt, \mbt^{\dag})\) is defined in (\ref{eq:covfn}) with 
\(f(\cdot; \mbl_0)\) replacing \(f(\cdot; \mbl)\).
\end{theorem}
{\em Proof.} The result follows immediately from the observations that
\begin{eqnarray*}
  n^{-1} \, \cov \left[ S(\mbt|\mbl_0), \, S(\mbt^{\dag}| \mbl_0)
  \right] &=& \mC(\mbt, \mbt^{\dag}), \\ n^{-1} \, \cov \left[\sum_{i
  = 1}^n \nabla \, l(\mbl_0|x_i), \sum_{i = 1}^n \nabla \,
  l(\mbl_0|x_i) \right] &=& {\bf I}(\mbl_0) \\ \mbox{and} \quad n^{-1}
  \, \cov\left[S(\mbt|\mbl_0), \, \sum_{i = 1}^n \nabla \,
  l(\mbl_0|x_i) \right] &=& \mC(\mbt|\mbl_0). \quad \hfill
  \rule{2mm}{2mm}
\end{eqnarray*} 
\assumption
\label{as:cstar}
Suppose \(\mC^{\star}(\mbt, \mbt^{\dag})\) is continuous and $0 <
\mC^{\star}(\mbt, \mbt) < \infty$ for all $\mbt \in \mbT$. 

\begin{theorem}
Under assumptions \ref{as:finitev} through \ref{as:cstar}, 
\begin{eqnarray*}
  \supt \frac{S(\mbt|\widehat{\mbl})}{\sqrt{n \, \mC^{\star}(\mbt,
  \mbt)}} \toinlaw \supt Z^{\star}(\mbt) \quad \mbox{as} \quad n \to
  \infty,
\end{eqnarray*}
where \(Z^{\star}(\mbt)\) is a Gaussian random field with the
covariance function 
\begin{eqnarray*}
  \rho^{\star}(\mbt, \mbt^{\dag}) := \frac{\mC^{\star}(\mbt,
  \mbt^{\dag})}{\sqrt{ \mC^{\star}(\mbt, \mbt) \,
  \mC^{\star}(\mbt^{\dag}, \mbt^{\dag})} }.
\end{eqnarray*}
\end{theorem}
{\em Proof.} First, the result holds for the process (\ref{eq:slth0})
(which is similar to Theorem \ref{th:maxsstar}). Next, the result
follows from Theorem \ref{th:nucov}. \hfill \rule{2mm}{2mm}  

We apply the results of Theorem \ref{th:sig3d} to the case of
one-dimensional \(\mbT\):
\begin{theorem} 
\label{th:nutail} 
The tail probability is expressed as $\mP\left( \supt Z(\mbt) \ge c
\right) = \kappa_0/(2 \pi) \; \mP(\chi^2_2 \geq c^2) + (\ell_0/4) \;
\mP(\chi^2_1 \geq c^2) + o[c^{-1} \exp(-c^2/2)]$ with 
\begin{eqnarray*}
  \k_0 = \int_{\mbT} \left. \left[ \frac{\partial^2}{\partial
  \mbt \, \partial \mbt^{\dag} } \rho^{\star}(\mbt, \mbt^{\dag})  
    \right]^{1/2} \right|_{\mbt^{\dag} = \mbt} \, 
    d \mbt 
\end{eqnarray*}
and \(\ell_0 = 2\).
\end{theorem}

The covariance function and \(\k_0\) depend on \(\mbl_0\); hence, 
cannot be evaluated directly. However, replacing \(\mbl_0\) by
\(\widehat{\mbl}\) yields a consistent estimator for $\mbl_0$. Just as
in the case of a fixed $\mbl$, the condition $\mC^{\star}(\mbt, \mbt) >
0$ for all $\mbt \in \mbT$ (part of assumption \ref{as:cstar}) will be
violated in the context of finite mixture models. However, one cannot
handle the singularities in a nice fashion and they are best treated
on a case-by-case basis. In particular, (1) there may be multiple
singularities, corresponding to each component of the mixture 
model under $\H_0$ and (2) in some cases the singularities lead to
discontinuities (as described earlier) while in other cases the
singularities are removable.

\setcounter{equation}{0} \def\theequation{\thesection.\arabic{equation}}

\section{Testing for the Order of a Mixture Model}
\label{sec-mixture}

In this section, building on the perturbation theory, we derive
results for the long-pending problem of testing for the order of a
mixture model while achieving the following goals: (1) Demonstrating
how the existing results for a special class of mixtures can be derived
from our general theory, (2) obtaining explicit and flexible
expressions for the geometric constants in the asymptotic tail
probability and (3) a careful examination of the singularities of the
score process that routinely occur in mixture models.

\subsection{Mixtures of Binomial Distributions}

Discrete mixtures for a random variable $X$ assuming a finite set of
values (e.g., \(0, \ldots, b\)) are of special interest, since the
data can be summarized by the bin counts \(N_0, \ldots, N_b\). The
loglikelihood and the score process \(S(\theta)\) depend on the data
only through these values. After appropriate centering and scaling, it
is easy to verify that the bin counts have an asymptotic $b$-variate
multivariate normal distribution. Consequently, the score process
$S(\theta)$ must have a finite Karhunen-Lo\`{e}ve expansion.

Consider the case of \(b = 2\) and a mixture of \(\mbox{Binomial}(2,
\theta)\) distributions with $\theta \in [0, 1]$. That is, our
interest is in testing $\H_0\!: \eta = 0$ against
$\H_1\!: \eta > 0$ and $\psi(x, \theta)$ is assumed to have a
$\mbox{Binomial}(2, \theta)$ distribution expressed as 
\begin{eqnarray*} 
  \psi(x; \theta) = \begin{cases} (1 - \theta)^2 & \mbox{if} \quad x = 0 \\ 2
  \theta \, (1 - \theta) & \mbox{if} \quad x = 1 \\ \theta^2 &
  \mbox{if} \quad x = 2 \end{cases} 
\end{eqnarray*} 
with the null density \(\psi(\cdot; \lambda)\) for some \(\lambda \in [0,
1]\).  Therefore, the perturbation model can be expressed as $p(x;
\eta, \lambda, \theta) = (1 - \eta) \; \psi(x, \lambda) +
\eta \; \psi(x, \theta)$.

\noindent {\em Case 1:} Assume $\lambda$ is known and $\theta$ is
unknown. The score process \begin{eqnarray*} S(\theta) = N_0 \,
\frac{(1 - \theta)^2}{(1 - \lambda)^2} + N_1 \, \frac{\theta \, (1 -
\theta)}{\lambda (1 - \lambda)} + N_2 \, \frac{\theta^2}{\lambda^2} -
n.  \end{eqnarray*} Since \(N_1 = (n - N_0 - N_2)\), the score process
reduces to \begin{eqnarray*} n^{-1/2} \, S(\theta) = Z_0 \, \frac{(1 -
\theta)(\lambda - \theta)}{(1 - \lambda)^2 \lambda} + Z_2 \,
\frac{\theta (\theta - \lambda)}{\lambda^2 \, (1 - \lambda)} =
{\mathfrak c}_0(\theta) \, Z_0 + {\mathfrak c}_2(\theta) \, Z_2,
\end{eqnarray*} where \(Z_0 = n^{-1/2} [N_0 - n \, (1 - \lambda^2)]\)
and \(Z_2 = n^{-1/2} (N_2 - n \, \lambda^2)\). The vector
\([{\mathfrak c}_0(\theta), \, {\mathfrak c}_2(\theta)]^T\) traces a
smooth curve through the origin at \(\theta = \lambda\). The
normalized score process $S^{\star}(\theta)$ has the flip property
discussed earlier.

The random variables \(Z_0\) and \(Z_2\) are correlated; hence,
explicit representation of \(S(\theta)\) in terms of the uncorrelated
random variables is quite messy. However, the corresponding manifold
\(\M\) consists of two arcs on the unit circle and the one-dimensional
volume of the tube is $\k_0 = \cos^{-1}(\mathfrak{r}_0) +
\cos^{-1}(\mathfrak{r}_1)$, where \(\mathfrak{r}_0 = \cor[S(0), -
S^{\prime}(\lambda)]\), \(\mathfrak{r}_1 = \cor[S(1),
S^{\prime}(\lambda)]\) and \(\ell_0 = 4\). Note that
\(\mathfrak{r}_0\) and \(\mathfrak{r}_1\) can be evaluated explicitly
based on $\mV(Z_0) = (1 - \lambda)^2 \, \lambda \, (2 - \lambda),
\cov(Z_0, Z_2) = - \lambda^2 \, (1 - \lambda)^2$ and $\mV(Z_2) =
\lambda^2 \, (1 - \lambda) \, (1 + \lambda)$. After some algebra, it
is easy to verify that \(\mathfrak{r}_0 = \sqrt{2 \lambda/(1 +
  \lambda)}\) and \(\mathfrak{r}_1 = \sqrt{2 (1 - \lambda)/(2 -
  \lambda)}\). Since $\M$ consists of two arcs on a unit circle, the
exact asymptotic null distribution of \(\mT\) is obtained using the
method of \cite{uusipaikka:83}.

\noindent {\em Case 2:} Assume that both $\lambda$ and $\theta$ are
unknown. Consider the MLE of \(\lambda\), \(\widehat{\lambda} = (N_1 +
2 N_2)/(2n) = (n + N_2 - N_0)/(2n)\), so that \(Z_0 = Z_2 = (N_0 +
N_2)/2 - n/4 - (N_2 - N_0)^2/(4n)\) and \(S(\theta|\widehat{\lambda})
= Z_0 (\theta - \widehat{\lambda})^2/[\widehat{\lambda}^2 \, (1 -
\widehat{\lambda})^2]\). In this case, the normalized score process is
constant and hence the manifold \(\M\) consists of a single
point. Therefore, \(\k_0 = 0\) and \(\ell_0 = 2\) resulting in a
distribution of \((0.5 \, \chi_0^2 + 0.5 \, \chi_1^2)\), where
$\chi_0^2$ is a degenerate distribution with all its mass at
zero. This is the special case derived by Lindsay (1995, p.\
95). \cite{shapiro:85} referred to this mixture of chi-square 
distributions with differing degrees of freedom as {\em chi-bar}
distribution.

\subsection{Mixtures of Exponential Family of Densities}
\label{sec-expf}

Suppose that \(\psi(x; \mbt)\) belongs to an exponential family of
densities so that $\psi(x; \mbt) = \exp[\mbt^T x - \varphi(\mbt)] \,
\psi_0(x)$. The null density is \(f(\cdot; \mbl)\) for some
\(\mbl\). 

\noindent {\em Case of Fixed $\mbl$:} The covariance function becomes
\begin{eqnarray*}
  \mC(\mbt, \mbt^{\dag}) &=& \int \exp[(\mbt + \mbt^{\dag} - \mbl)^T
  x + \varphi(\mbl) - \varphi(\mbt) - \varphi(\mbt^{\dag})] \,
  \psi_0(x) \, dx - 1 \\ &=& \exp[\varphi(\mbt + \mbt^{\dag} -
  \mbl) + \varphi(\mbl) - \varphi(\mbt) - \varphi(\mbt^{\dag})] -
  1.
\end{eqnarray*}
If $\psi(\cdot; \mbt)$ has a multivariate normal distribution with a
mean vector $\mbt$ and an identity variance covariance matrix, it
follows that \( \varphi(\mbt) = \|\mbt\|^2/2\) and 
\begin{equation}
  \mC(\mbt, \mbt^{\dag}) = \exp[\ip{\mbt - \mbl}{\mbt^{\dag} -
  \mbl}] - 1. 
  \label{eq:normcov}
\end{equation}

%Only for this special case, the above covariance function was obtained
%by \cite{chen1:01}. 
Consider the special case of $d = 1$.  The critical values are
obtained using Theorem \ref{th:sig3d} and the one-dimensional volume
of $\M$ has the following explicit expression when $\lambda = 0$:
 \begin{eqnarray*} \k_0 =
\int_{\Theta} \frac{ \left[ \exp(2 \theta^2) - (1 + \theta^2) \,
    \exp(\theta^2) \right]^{1/2}}{ [\exp(\theta^2) - 1] } \, d \theta.
\end{eqnarray*} 
The normalized score process again has the flip property
(\ref{eq:flip}) and \(\ell_0 = 4\).

\noindent {\em Case of Unknown $\mbl$:} Straightforward calculations
show that the covariance function (\ref{eq:covfn2}) in Theorem
\ref{th:nucov} becomes 
\begin{eqnarray*}
  \mC^{\star}(\mbt, \mbt^{\dag}) &=&
    \mC(\mbt, \mbt^{\dag}) - [\varphi^{\prime}(\mbt) -
    \varphi^{\prime}(\mbl)]^T \, \left[ \varphi^{\prime
    \prime}(\mbl) \right]^{-1} \, [\varphi^{\prime}(\mbt^{\dag}) -
    \varphi^{\prime}(\mbl)], 
\end{eqnarray*}
since $\mC(\mbt|\mbl) = \mE_{\mbt}[X - \varphi^{\prime}(\mbl)] =
\varphi^{\prime}(\mbt) - \varphi^{\prime}(\mbl)$ and ${\bf I}(\mbl) =
\varphi^{\prime \prime}(\mbl)$. 

In the case of a univariate normal distribution, the volume of the
one-dimensional manifold becomes
\begin{eqnarray*} 
  \k_0 = \int_{\Theta} \frac{ \left[ \exp\{2 (\theta - \lambda)^2\} +
      1 - \exp\{(\theta - \lambda)^2\} \, \left\{ 2 + (\theta -
        \lambda)^4 \right\} \right]^{1/2}}{\left[ \exp\{(\theta -
      \lambda)^2 \} - 1 - (\theta - \lambda)^2 \right] } \, d \theta.
\end{eqnarray*} 
The normalized score process has a singularity at \(\theta =
\widehat{\theta}\); however, the precise behavior at this point needs
careful consideration, which is presented next. In the neighborhood of
\(\widehat{\theta}\), we have
\begin{equation}
\label{eq:quadscore}
  S(\theta) = S(\widehat{\theta}) + (\theta - \widehat{\theta}) \,
  S^{\prime}(\widehat{\theta}) + \frac{1}{2} \, (\theta -
  \widehat{\theta})^2 \, S^{\prime \prime}(\widehat{\theta}) + 
  o\left[(\theta - \widehat{\theta})^2 \right].
\end{equation}
Note that \(S(\widehat{\theta}) = S^{\prime}(\widehat{\theta}) = 0\)
(since the latter is simply the score equation defining
\(\widehat{\theta}\)). By continuity, \(S^{\prime
\prime}(\widehat{\theta}) = \varphi^{\prime \prime}(\lambda) + o(1)\);
hence, the normalized score process becomes $S^{\prime
\prime}(\lambda)/\sqrt{\mV[S^{\prime \prime}(\lambda)]} + o(1)$ in
the neighborhood of \(\lambda\). This is continuous so there is
no flip at \(\theta = \widehat{\theta}\). The manifold \(\M\) for this
process is a single segment and \(\ell_0 = 2\).

\subsection{Testing for \(m\) versus $(m + q)$ Component Mixture Model}

One of the important applications of the perturbation theory is in
building finite mixture models formed from a broad class of smooth
densities. First, consider testing \begin{eqnarray*} \H_0\!: 
\mbox{$m$-component mixture} \quad \mbox{against} \quad \H_1\!: \mbox{$(m +
  1)$-component mixture} \end{eqnarray*} when mixtures are formed from
  any smooth families, including discrete, continuous and multivariate
  densities. Under the mixture model framework, the null model
\(f(\cdot; \mbl)\) is the \(m\) component mixture ${\rm g}(x; \Q_m) =
\sum_{j = 1}^m \beta_j \, \psi(x; \mbt_j)$, where $\Q_m = (\mbt^T,
\mbb^T)^T$ while the alternative is the \((m + 1)\)st component. We
consider two cases: (1) The support point vectors $\mbt$s are fixed
and only the mixing weight vector $\mbb$ is estimated and (2) 
$\mbt$s and $\mbb$ are estimated.

\noindent{\em Case 1:} Assume $\mbt$ is fixed and the goal is to
estimate $\mbb$. The likelihood surface is concave in \(\mbb\) and the MLEs
satisfy
\begin{equation}
  \label{eq:betaj}
  S(\mbt_j|\widehat{\mbb}) = \sum_{i = 1}^n \frac{ \psi(x_i; \mbt_j)
  }{ {\rm g}(x_i; \widehat{\Q}) } - n = 0 \quad \mbox{for all} \quad j
  = 1, \ldots, m   
\end{equation}
provided that the solution satisfies \(0 < \widehat{\beta}_j < 1\)
(otherwise, some components are set to zero). The MLE satisfies the
conditions of Section \ref{sec-nuisance}, provided that \(\beta_j > 0
\) for each \(j\). The covariance function is determined based on the
result in Theorem \ref{th:nucov}. 

The set of equations in (\ref{eq:betaj}) implies that the normalized
score process has a singularity at each \(\mbt_j\). Using an argument
similar to (\ref{eq:flip}), the process flips at each of these
points. 

\noindent{\em Case 2:} The goal is to estimate both $\mbt$ and
$\mbb$. Note that each support point is of dimension $d$. The
equations defining the MLEs become 
\begin{eqnarray} 
   \label{eq:scd}  
   S(\mbt_j| \wQ_m)\Big|_{\mbt_j = \wmbt_j} &=& 0 
   \quad \mbox{and} \quad S^{\prime}(\mbt_j| \wQ_m)\Big|_{\mbt_j = \wmbt_j}
   = {\bf 0}
\end{eqnarray}
for all $j = 1, \ldots, m$. Note that for \(d > 1\), the above
equation is a vector. Using an expansion similar to
(\ref{eq:quadscore}), around each of the true support points, it is
easy to verify that all the singularities in the normalized score
process are removable.

Consistent estimators of the nuisance parameters are required to apply
the results of Section \ref{sec-nuisance}. This is achieved by
imposing an order constraint on the support point vectors \(\mbt_j\)
and a corresponding constraint on the estimators. Under these
constraints, the approximate critical values are obtained from
Theorems \ref{th:nucov} and \ref{th:nutail}. 

\begin{figure}[htbp]
  \centerline{\scalebox{0.6}{\includegraphics{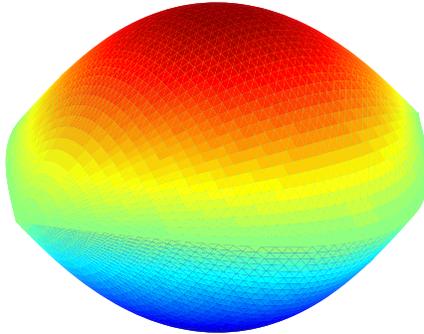}}}
\caption{Manifold for testing \(m\) versus $(m + 2)$ components in
mixture models. The manifold has two corners, two edges and two
boundary faces.} 
\label{fig:wedge}
\end{figure}

\vspace{0.1in}
\noindent {\em General case:} Consider the more general problem of
testing $\H_0\!: \mbox{$m$-component mixture}$ against $\H_1\!:
\mbox{$(m + q)$-component mixture for $q = 1, 2, \ldots$}$. For this
case, Theorem \ref{th:sig3d} is still applicable and the score process
is easy to derive (see \citep{pilla2:03} for details). Suppose \(d =
1\) and \(\mbT\) is an interval, then the manifold has two corner points and 
two edges with two boundary faces as shown in Fig.\ \ref{fig:wedge}. 

\subsection{Mixtures of Bivariate Normal Distributions}
\label{sec-bivariate}

In this section, we consider the bivariate mixture testing problem so
that $d = 2, {\bf x} = (x_1, x_2)^T$ and \(\mbt = (\theta_1,
\theta_2)^T\). To the best of the authors' knowledge, this is the
first attempt at testing for mixtures of multivariate
distributions. Assume $f(\cdot; \mbl)$ is a bivariate standard normal
density and $\psi(\cdot; \mbt)$ is a bivariate normal density with
mean $\mbt$ and an identity covariance matrix. From equation
(\ref{eq:normcov}), it is easy to verify that the covariance function
can be explicitly expressed as \( \mC(\mbt, \mbt^{\dag}) =
\exp[\ip{\mbt}{\mbt^{\dag}}] - 1 \). Suppose \(\mbT\) is a disk of
radius \(\varrho_1 > 0\), so that \begin{eqnarray*} \mT = \sup_{0 < 
\|\mbt\| \le \varrho_1} \frac{S(\mbt)}{\sqrt{n \, \mC(\mbt, \mbt)} }.
  \end{eqnarray*} In order to address the singularity at \(\|\mbt\| =
  0\), first consider the supremum over \(\varrho_0 \le \|\mbt\| \le
  \varrho_1\), where \(0 < \varrho_0 < \varrho_1\) and next let
  \(\varrho_0 \to 0\). Under the polar coordinate parameterization of
  \(\mbt = [\varrho \cos(\omega),
\varrho \sin(\omega)]^T \), with the covariance function expressed as
\(\mC(\mbt, \mbt^{\dag}) = \exp[\varrho \, \varrho^{\dag} \cos(\omega
- \omega^{\dag})] - 1\), it follows that
\begin{eqnarray*}
  \k_0 &=& \int_{\varrho_0}^{\varrho_1} \int_0^{2\pi} [\mC(\mbt, \mbt)
  ]^{-3/2} \, \det \begin{pmatrix} \exp(\varrho^2) - 1 & \varrho \,
  \exp(\varrho^2) & 0 \\ \varrho \, \exp(\varrho^2) & (1 + \varrho^2)
  \, \exp(\varrho^2) & 0 \\ 0 & 0 & \varrho^2 \, \exp(\varrho^2)
  \end{pmatrix}^{1/2} \, d \, \omega \, d \varrho \\ &=& 2 \pi
  \int_{\varrho_0}^{\varrho_1} \left[ \frac{\varrho^2 \, \exp(3
  \varrho^2) - \varrho^2 \, (1 + \varrho^2) \, \exp(2
  \varrho^2)}{\{\exp(\varrho^2) - 1\}^3 } \right]^{1/2} \, d \varrho.
\end{eqnarray*}
The integrand has a finite limit as \(\varrho \to 0\); therefore, the
integral is still valid when \(\varrho_0 = 0\).  

Next, we consider the boundaries at \(\varrho = \varrho_0\) and
\(\varrho = \varrho_1\). For an arbitrary \(\varrho\), the length of the
boundary is 
\begin{eqnarray*}
  \ell_0 = \int_0^{2\pi} [\mC(\mbt, \mbt)]^{-1} \, 
    \det \begin{pmatrix}
    \exp(\varrho^2) - 1 & 0 \\ 0 & \varrho^2 \, \exp(\varrho^2)
    \end{pmatrix}^{1/2} \, d \omega
  = 2 \pi \sqrt{\frac{ \varrho^2 \, \exp(\varrho^2)
  }{ [\exp(\varrho^2) - 1]} }.
\end{eqnarray*}
Therefore, 
\[
  \ell_0 = 
   2 \pi \left[ \sqrt{ \frac{ \varrho_0^2 \, \exp(\varrho_0^2)
   }{[\exp(\varrho_0^2) - 1]} } + \sqrt{ \frac{
   \varrho_1^2 \, \exp(\varrho_1^2) }{ 
    [\exp(\varrho_1^2) - 1] } } \right]
    \longrightarrow 2 \pi \left[ 1 + \sqrt{ \frac{ \varrho_1^2 \,
    \exp(\varrho_1^2) }{ \{\exp(\varrho_1^2)- 1\} } } \right] 
\]
as $\varrho_0 \to 0$. The contribution from the inner boundary does
not disappear as \(\varrho_0 \to 0\), instead it converges to \(2
\pi\). This implies that the manifold \(\M\) corresponding to this
process has a hole and \(\M\) has an Euler-Poincare characteristic of
\(\mathcal{E} = 0\). The tail-probability approximation of Theorem
\ref{th:sig3d} simplifies to 
 \begin{eqnarray*} 
   \mP \left(\supt Z(\mbt) \geq c \right) &\approx& \frac{\k_0}{4 \pi}
   \mP \left( \chi_3^2 \ge c^2\right) + \frac{\ell_0}{4 \pi} \mP
   \left(\chi_2^2 \ge c^2 \right) - \frac{\k_0}{4 \pi} \mP \left(
     \chi_1^2 \ge c^2 \right) \\
   &=& \frac{\k_0}{2 \sqrt{2\pi}} \, c \, \exp(-c^2/2) +
   \frac{\ell_0}{4 \pi} \exp(-c^2/2) \quad \mbox{as} \quad n \to
   \infty.  
\end{eqnarray*}

The interior hole occurs in any two-parameter problem, as the next
lemma demonstrates.  

\begin{figure}[htbp]
\centerline{\scalebox{0.6}{\includegraphics{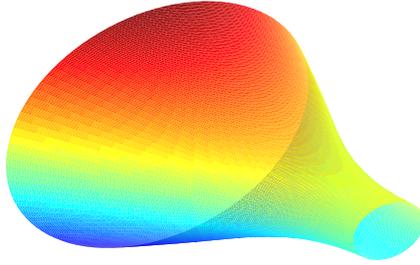}}}
\caption{Manifold for the bivariate normal mixture testing
problem. The cylindrical manifold has two boundaries: a circle with
circumference \(2 \pi\), corresponding to \(\varrho = 0\), and a
larger (high dimensional) ring corresponding to \(\varrho =
\varrho_1\).} 
\end{figure}

\begin{lemma} 
\label{th:bound}
Suppose \(\mbt\) is of dimension \(d = 2\) and there exists a $\mbl$
such that $f(\cdot; \mbl) = \psi(\cdot; \mbt)$. The normalized
score process \(S^{\star}(\mbt)\) has a singularity at $\mbt = \mbt_0$
and correspondingly, the manifold \(\M\) has a hole. The length of the
interior boundary of $\M$ is $2 \pi$. 
\end{lemma}
{\em Proof.} A Taylor series expansion yields
\begin{eqnarray*}
  S(\mbt) = \left < \mbt - \mbt_0, \, S^{\prime}(\mbt_0)
  \right > + o(\|\mbt - \mbt_0\|) \quad \mbox{as} \quad \mbt \to
  \mbt_0. 
\end{eqnarray*}
Let \(\matr{R}\) be a matrix such that \(\cov[S^{\prime}(\mbt_0)] = n
\, \matr{R}^T \, \matr{R}\). Then the normalized score process becomes
\begin{eqnarray*}
  S^{\star}(\mbt) = \frac{S(\mbt)}{\sqrt{n \, \mC(\mbt, \mbt)}}
     = \left < \frac{\matr{R}(\mbt - \mbt_0) }{
     \|\matr{R}(\mbt - \mbt_0)\| }, \frac{\matr{R}^{-1}
     S^{\prime}(\mbt_0)}{\sqrt{n}} 
     \right > + o(\|\mbt - \mbt_0\|).
\end{eqnarray*}
As \(\mbt\) varies in a small circle around \(\mbt_0\),
the boundary of the manifold $\M$, \( \matr{R}(\mbt -
\mbt_0)/ \|\matr{R}(\mbt - \mbt_0)\| \), becomes the unit
circle in \(\R^2\) which has length \(2 \pi\).  \hfill
\rule{2mm}{2mm}    

For \(d = 1 \), the manifold \(\M\) has \((m + 1)\) segments so that
\(\ell_0 = 2(m + 1)\). Approximate critical values are obtained based
on Theorem \ref{th:sig3d} and \(\k_0\) is evaluated using numerical
integration. For \(d = 2\), the manifold \(\M\) has \(m\) holes with
each hole contributing \(2 \pi\) to the total length of the boundary
\(\ell_0\). The Euler-Poincare characteristic of \(\M\) is therefore
$(1 - m)$. For the result in Theorem \ref{th:sig3d}, the constant
\(\k_0\) and the length of the outer boundary are found using a
bivariate and univariate numerical integrations, respectively.

\subsection{Simulation Experiments}

In order to demonstrate the power of the proposed methods, we present
two simulation studies and illustrate the process of building mixture
models. 

\begin{figure}[htb]
\psfrag{Normalized Score}{\hspace{0.5in}$S^{\star}(\theta)$}
\psfrag{lambda}{\large $\theta$}
\centerline{\scalebox{0.6}{\includegraphics{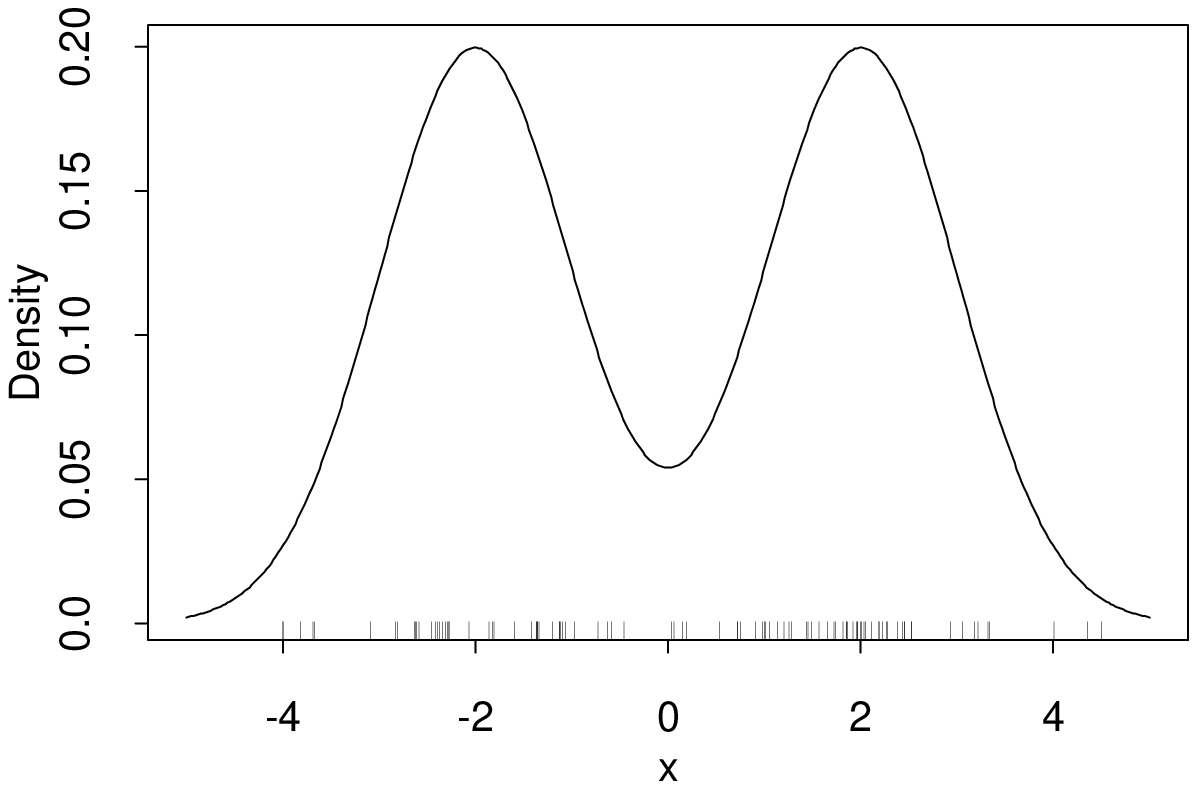} \hfill
\includegraphics{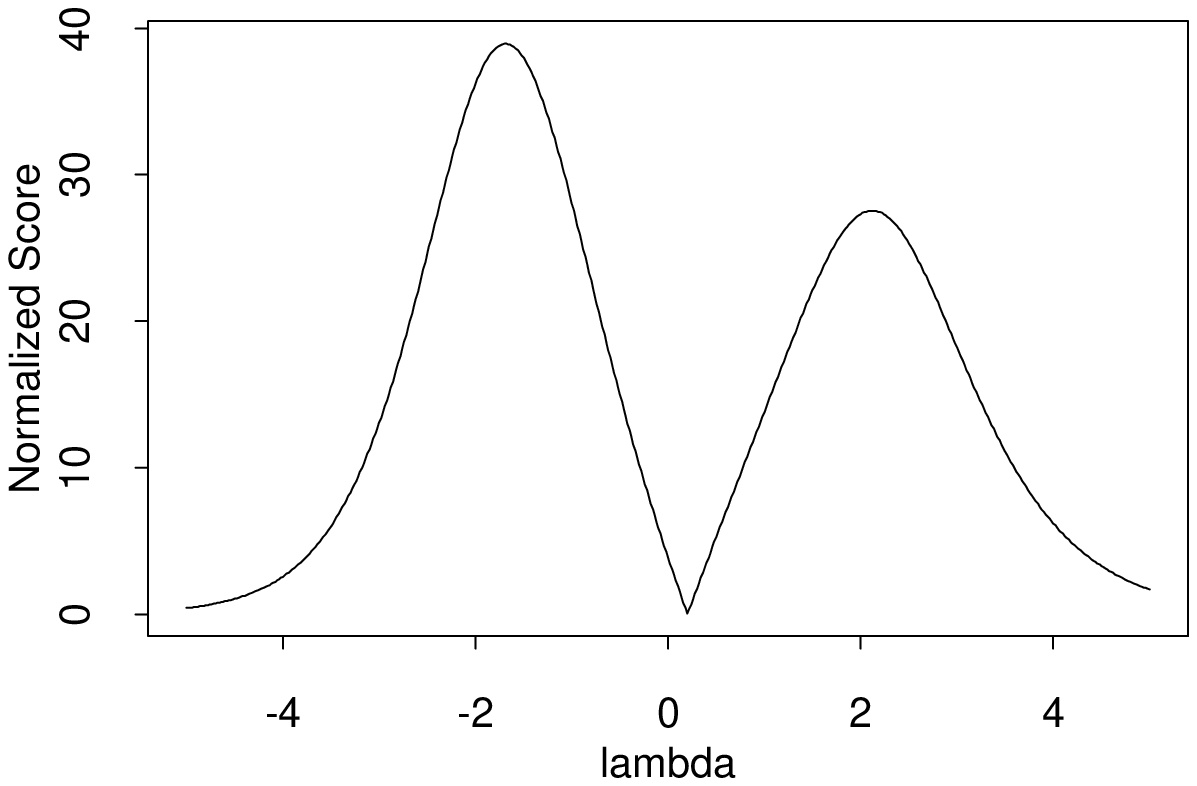}} }
\centerline{\scalebox{0.6}{\includegraphics{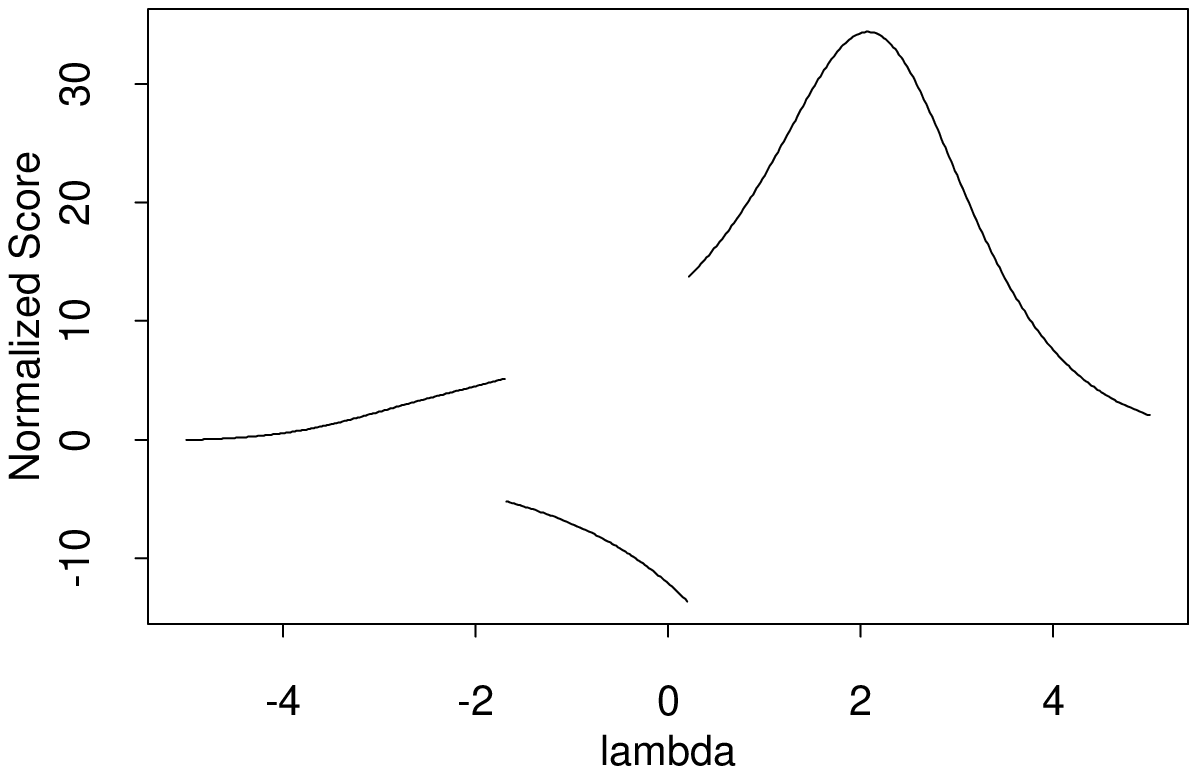} \hfill
\includegraphics{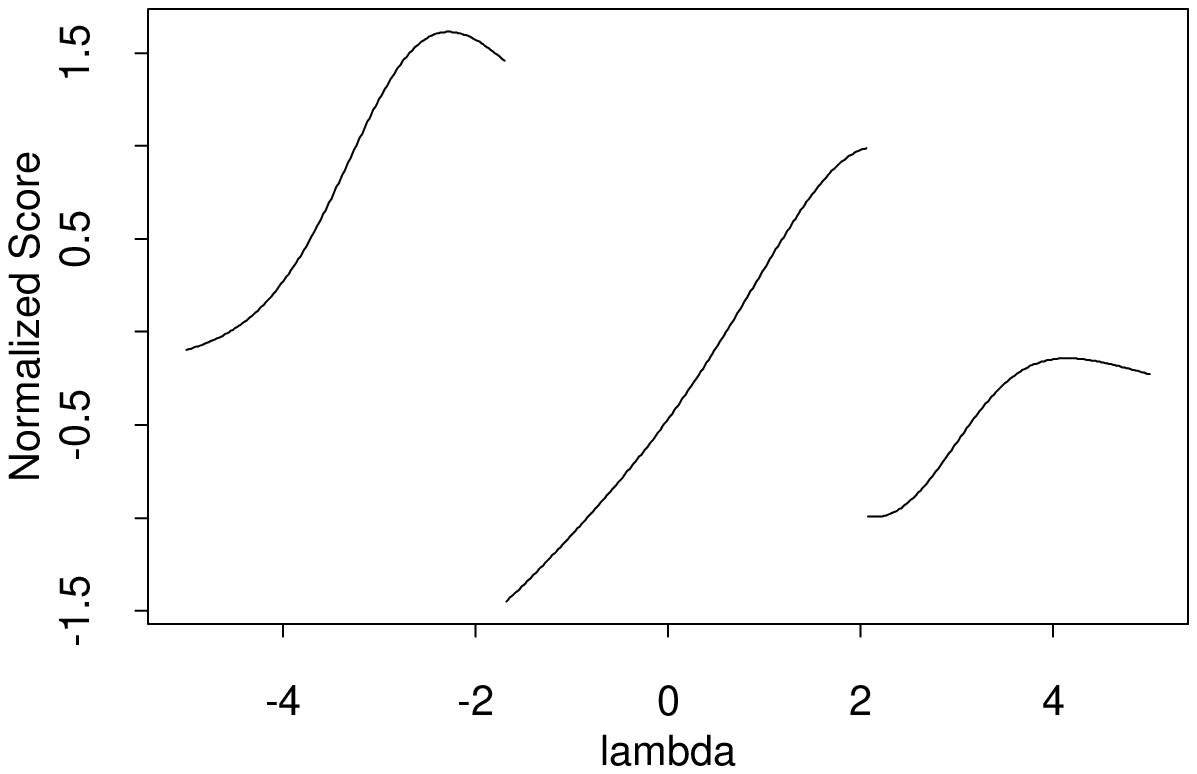}} }
\label{fig:data}
\caption{(a) Simulated data from the two-component normal mixture: \(0.5
  N(-2, 1) + 0.5 N(2, 1)\); (b) $S^{\star}(\theta)$ for the
  one-component mixture model; (c) $S^{\star}(\theta)$ for the
  two-component mixture model; (d) $S^{\star}(\theta)$ for the model
  with a third component included and the first component removed.}
\end{figure}
We consider the simulated dataset shown in Fig.\ \ref{fig:data}(a),
consisting of a sample of size $n = 100$ drawn from the two-component
normal mixture model $0.5 N(-2, 1) + 0.5 N(2, 1)$. The model building
process starts with the first component at the sample mean \(
\widehat{\theta}_1 = \overline{X} = 0.20322\). The starting model is
obviously a poor fit for the dataset. Fig.\ \ref{fig:data}(b) presents
the fitted normalized score process $S^{\star}(\theta)$, showing two
peaks in the vicinity of the true mixture components.  An application
of the volume-of-tube formula in (\ref{eq:tubechis}) to this model
yields \(\k_0 = 5.72\) and \(\ell_0 = 4\) with the critical value of
\(c = 2.518\) at the 5\% level. Clearly the peaks are highly
significant. A second component at \(\theta_2 = -1.68929\) [the
location of the larger left peak in Fig.\ \ref{fig:data}(b)] is
included in the model and the vector of estimated mixing weights is 
\(\widehat{\mbb} = (0.67315, 0.32685)^T\). The incorrect first
component \(\theta_1\) still dominates the fitted mixture model.

The normalized score process relative to the two-component mixture is
shown in Fig.\ \ref{fig:data}(c). The striking feature of this plot is
the two discontinuities at the fitted components \(\theta_1 = 0.203\)
and \(\theta_2 = -1.689\). These discontinuities occur due to the
zeroes of the covariance function since \(\mC^{\dag}(\theta_1,
\theta_1) = \mC^{\dag}(\theta_2, \theta_2) = 0\) which in turn
corresponds to the singularities in $S^{\star}(\theta)$.  The manifold
\(\M\) for this process has three pieces so that $\k_0 = 5.082$ and
$\ell_0 = 6$. The critical value \(c = 2.571\) and the right peak is
still highly significant. The maximum occurs at \(\widehat{\theta}_3 =
2.07328\) which is included as a third component in the model. Since
\(\widehat{\mbb} = (0, 0.45616, 0.54384)^T\), the first component is
removed from the model. For the two-component mixture model with
\(\widehat{\theta}_2\) and \(\widehat{\theta}_3\), the constants
\(\k_0 = 5.082\) and \(\ell_0 = 6\) yielding \(c = 2.571\). Fig.\
\ref{fig:data}(d) presents the process $S^{\star}(\theta)$ and it is
entirely below the critical value \(c\); therefore, the two-component
mixture model with $\widehat{\theta}_2$ and $\widehat{\theta}_3$ is
the final fitted model.

The true density is chosen as $p(x; \eta, \theta) = 0.5 (1 - \eta) \,
\psi(x; -2) + \eta \, \psi(x; 0) + 0.5 (1 - \eta) \, \psi(x; 2)$ for
\(\eta \in \{0, 0.1, 0.2\}\) and \(\psi(\cdot; \theta)\) is the normal
density with mean \(\theta\) and unit variance. This density has two
large with well separated components and our goal is to test for the
presence of the poorly separated third component. We present
simulation studies using 1000 data sets under the following three
different scenarios: \\ 
{\em Model 1:} \(f(x; \mbl) \equiv {\rm g}(x; \Q) = [0.5
\, \psi(x; -2) + 0.5 \, \psi(x; 2)]\) is completely specified. \\
{\em Model 2:} \(f(x; \mbl) \equiv {\rm g}(x; \Q) =
[\beta_1 \, \psi(x; -2) + \beta_2 \, \psi(x; 2)] \), where $\beta_1$
and $\beta_2$ are estimated. \\ 
{\em Model 3:} \(f(x; \mbl) \equiv {\rm g}(x; \Q) = [\beta_1 \,
\psi(x; \theta_1) + \beta_2 \, \psi(x; \theta_2)] \), where $\beta$s
and $\theta$s are estimated. 

\begin{table}[htb]
\vspace{-0.1in}
\caption{Rejection rates for three different null models under three 
different perturbation sizes based on 1000 simulation studies.} 
\label{tab:simul}
\vspace{.5cm} 
\centerline{
\begin{tabular}{|c|ccc|ccc|} \hline 
   & \multicolumn{3}{c|}{$n = 200$} & \multicolumn{3}{c|}{$n =
   1000$} \\ 
   Model & $\eta = 0.0$ & $\eta = 0.1$ & $\eta = 0.2$ & $\eta = 0.0$
	 & $\eta = 0.05$ & $\eta = 0.1$ \\ \hline
   1 & 79 & 537 & 975 & 74 & 636 & 990 \\
   2 & 78 & 583 & 985 & 76 & 673 & 992 \\
   3 & 74 & 292 & 588 & 61 & 371 & 817 \\ \hline 
\end{tabular} }
\end{table}
Table \ref{tab:simul} presents the rejection rates for 1000
simulations under two sample sizes. When \(\eta = 0\), $\H_0$ is true
and hence we expect the rejection rate to be close to the nominal
significance level of $5\%$. As \(\eta\) increases, the power
increases as expected. As the null assumptions are relaxed, the power
decreases which again is to be expected. The poor separation between
the components makes it difficult for the test to detect the third
component which is more prominent for model 3. Naturally, estimating the
nuisance parameters under the null model has an effect on the power of
the test.

\section{Discussion}
\label{sec-discuss}

In this article, we introduced a general class of models, perturbation
models, and proposed a test statistic (asymptotically equivalent to
the LRT statistic) based on the score process to detect the presence
of perturbation.  We derived general inferential theory for the
asymptotic null distribution of the test statistic for a class of
non-regular problems using the Hotelling-Weyl-Naiman volume-of-tube
formula. The resulting theory is extended to solve the long-pending
fundamental problem of testing for the mixture complexity, including
the case when the null model includes a set of nuisance
parameters. Our theory is applicable to a general family of mixture
models including the multivariate family of mixtures. Other
applications to the general theory include spatial scan analysis,
latent class models (employed in social research) and Rasch models
(employed in educational testing and survey sampling). The inferential
theory developed in this article provides a solution to an important
class of statistical problems involving loss of identifiability and/or
when some of the parameters are on the boundary of the parametric
space.

The explicit determination of the geometric constants appearing in the
tube formula are carried out using the \texttt{Libtube} software
\citep{loader:05}. Our theory is general enough to be applicable to
scalar or vector $\mbl$ and univariate or multivariate data. The
advantage of our approach is that the tube formula provides an elegant
approximation to the asymptotic null distribution compared to those
based on simulations or bootstrap based procedures.

\setcounter{equation}{0}
\def\theequation{\thesection.\arabic{equation}}

\section{Proofs}
\label{sec-proof}

In this section we provide proofs of the main theorems.  As before,
notation $\prime$ is used to denote derivative with respect to the
appropriate term. 

\noindent {\em Proof of Theorem \ref{th:lrtscore}.} Let 
\begin{eqnarray*} 
  K(\eta, \mbt) = \sum_{i = 1}^n \log \left[1 + \frac{\eta \, \left\{\psi(x_i;
      \mbt) - f(x_i; \mbl) \right\}}{f(x_i; \mbl)} \right].  
\end{eqnarray*} 
The LRT statistic becomes $\sup_{\mbt, \eta > 0} K(\eta, \mbt)$. For
any \(\eta > 0\), a Taylor series expansion yields 
\begin{eqnarray*}
  K(\eta/\sqrt{n}, \mbt) &=& K(0, \mbt) + \frac{\eta}{\sqrt{n}}
  K^{\prime}(0, \mbt) + \frac{\eta^2}{2 n} K^{\prime
    \prime}(\eta^{\star}, \mbt) \quad \mbox{for} \quad 0 \leq
  \eta^{\star} \leq \frac{\eta}{\sqrt{n}} \\ 
&=& \frac{\eta}{\sqrt{n}} \sum_{i = 1}^n \left[ \frac{\psi(x_i;
    \mbt)}{f(x_i; \mbl)} - 1 \right] \\
&& - \, \frac{\eta^2}{2 n} \sum_{i = 1}^n \left[ \frac{ \left\{
      \psi(x_i; \mbt) - f(x_i; \mbl) \right\}^2/ \{f(x_i; \mbl)\}^2}{1
    + \eta^{\star} \{\psi(x_i; \mbt)/f(x_i; \mbl) - 1\} } \right].
\end{eqnarray*} 
Under an implicit assumption that convergence statements are uniform
in $\mbt$ for bounded sets and from the results in \cite{rubin:56}, it
follows that
\begin{eqnarray*} 
   K^{\prime \prime}(\eta^{\star}, \mbt) = \sum_{i = 1}^n \left[
     \frac{ \left\{ \psi(x_i; \mbt) - f(x_i; \mbl) \right\}^2/
       \{f(x_i; \mbl)\}^2}{1 + \eta^{\star} \{\psi(x_i; \mbt)/f(x_i;
       \mbl) - 1\} } \right] 
\end{eqnarray*} 
is uniformly converging to $\mC(\mbt, \mbt)$. Therefore, 
\begin{eqnarray*}
   K(\eta/\sqrt{n}, \mbt) &=& \frac{\eta}{\sqrt{n}} S(\mbt) -
  \frac{\eta^2}{2} \mC(\mbt, \mbt) + o_p(1), 
\end{eqnarray*} 
where the $o_p(1)$ term is uniform in $\eta$ and $\mbt$ on compact
sets. In effect, $\sup_{\eta \geq 0} K(\eta/\sqrt{n}, \mbt) =
(1/2) \max\{0, S^{\star}(\mbt)\}^2 + o_p(1)$. \hfill \rule{2mm}{2mm}

On the way to proving Theorem \ref{th:maxsstar}, we derive a series of
technical results.

\begin{lemma}
\label{lem:linb}
Let \(a(\theta)\) be a continuously differentiable function
on an interval \(\Theta\). Let \(a_{\star} = [a(\theta_1) -
a(\theta_0)]\). Then   
\begin{eqnarray*}
  \int_{\theta_0}^{\theta_1} [a^{\prime}(\theta)]^2 \, d \theta
    \ge \frac{a_{\star}^2}{|\theta_1 - \theta_0|},
\end{eqnarray*}
where $a^{\prime}(\theta) = d a(\theta)/d \theta$. 
\end{lemma}
{\em Proof.} Let \(\theta_{\star} = (\theta_1 - \theta_0) \) so that
\begin{eqnarray*}
  \int_{\theta_0}^{\theta_1} [a^{\prime}(\theta)]^2 \, d
  \theta &=& \int_{\theta_0}^{\theta_1} \left( a^{\prime}(\theta)
    - \frac{a_{\star}}{\theta_{\star}} + \frac{a_{\star}}{\theta_{\star}}
    \right)^2 \, d \theta \\ 
  &=& \int_{\theta_0}^{\theta_1} \left( a^{\prime}(\theta) -
  \frac{a_{\star}}{\theta_{\star}} \right)^2 \, d \theta +
  \int_{\theta_0}^{\theta_1} \frac{a_{\star}^2}{\theta_{\star}^2} \, d
  \theta + \frac{2 a_{\star}}{\theta_{\star}}\int_{\theta_0}^{\theta_1} \left(
  a^{\prime}(\theta) - \frac{a_{\star}}{\theta_{\star}} \right) \, d \theta.
\end{eqnarray*}
Note that the first integral is non-negative and the third one is
zero. \hfill \rule{2mm}{2mm} 

\begin{lemma}
\label{lem:blin2}
Suppose \(\theta_0 < \theta_2\) and \(a(\theta_0) = a(\theta_2) = 0\), then 
\begin{eqnarray*}
  \int_{\Theta} [a^{\prime}(\theta)]^2 \, d \theta
    \ge  \frac{4}{|\theta_2 - \theta_0|}
      \left(\sup_{\theta_0 \le \theta \le \theta_2} |a(\theta)|
      \right)^2.
\end{eqnarray*}
\end{lemma}
{\em Proof.} Suppose the supremum occurs at \((\theta_1, a_{\star})\) with
\(\theta_0 < \theta_1 < \theta_2\). An application of Lemma \ref{lem:linb}
separately over \([\theta_0, \theta_1]\) and \([\theta_1, \theta_2]\) yields  
\begin{eqnarray*}
  \int_{\Theta} [a^{\prime}(\theta)]^2 \, d \theta \ge
  \int_{\theta_0}^{\theta_2} [a^{\prime}(\theta)]^2 \, d \theta
    \ge a_{\star}^2 \left[ \frac{1}{(\theta_1 - \theta_0)} +
       \frac{1}{(\theta_2 - \theta_1)} \right]
    \ge \frac{ 4 a_{\star}^2 }{ (\theta_2 - \theta_0) }.
\end{eqnarray*} \hfill \rule{2mm}{2mm}

\begin{lemma}
\label{lem:linb3}
Suppose \(b(\theta)\) is continuously differentiable. For
\(\delta > 0\), let \(b_{\delta}(\theta)\) be the linear interpolant
between the points \(0, \pm \delta, \pm 2 \delta, \ldots\). Then
\begin{eqnarray*}
  \supt \Big| b_{\delta}(\theta) - b(\theta) \Big|^2
     \le \delta \int_{\Theta} [b^{\prime}(\theta)]^2 \, d \theta.
\end{eqnarray*}
\end{lemma}
{\em Proof.} Once again, let \(a_{\star}\) be the supremum. An
application of Lemma \ref{lem:blin2} to \(a(\theta) = [b_{\delta}(\theta)
- b(\theta)] \) yields 
\begin{eqnarray*}
  a_{\star}^2 &\le& \frac{\delta}{4} \int_{\Theta} \left[
  b^{\prime}_{\delta}(\theta) - b^{\prime}(\theta) \right]^2 \, d
  \theta \le \frac{\delta}{2} \int_{\Theta}
    \left[ \{b^{\prime}_{\delta}(\theta)\}^2 +
    \{b^{\prime}(\theta)\}^2 \right] \, d \theta \\ 
  &\le& \delta \int_{\Theta} [b^{\prime}(\theta)]^2 \, d \theta.
\end{eqnarray*}
The final inequality holds since \(\int_{\Theta}
[b^{\prime}_{\delta}(\theta)]^2 \, d \theta \le \int_{\Theta}
[b^{\prime}(\theta)]^2 \); this follows from the application of Lemma
\ref{lem:linb} between each pair of knots of
\(b_{\delta}(\cdot)\). \hfill \rule{2mm}{2mm}  

\begin{lemma}\label{lem:unibd}
Let \(Y(\mbt)\) be a stochastic process with continuously
differentiable sample paths and let \(Y_{\delta}(\mbt)\) be its
linear interpolant between points \(0, \pm \delta, \ldots\). 
Then
\[
  \mP \left( \supt |Y_{\delta}(\mbt) - Y(\mbt)| \ge
  \epsilon \right) \le \frac{\delta}{\epsilon^2} \, \mE
  \int_{\mbT} [Y^{\prime}(\mbt)]^2 \, d \mbt. 
\]
Uniform convergence holds if the expectation is finite:
\begin{equation}
  \label{eq:unifc}
  \lim_{\delta \to 0}
    \mP \left( \supt |Y_{\delta}(\mbt) - Y(\mbt)| \ge \epsilon
    \right) = 0 \quad \mbox{for all} \quad \epsilon > 0.
\end{equation}
\end{lemma}
{\em Proof.} From Lemma \ref{lem:linb3}, it follows that
\begin{eqnarray*}
  \mP \left( \supt |Y_{\delta}(\mbt) - Y(\mbt)| \ge \epsilon
  \right) &\le& \mP \left( \delta \int_{\mbT}
  [Y^{\prime}(\mbt)]^2 \, d \mbt \ge \epsilon^2 \right) \\ 
    &\le& \frac{\delta}{\epsilon^2} \, \mE \int_{\mbT}
    [Y^{\prime}(\mbt) ]^2 \, d \mbt ,
\end{eqnarray*}
where the last line follows from the Markov's inequality for any
non-negative random variable.  \hfill \rule{2mm}{2mm}

\begin{lemma}\label{lem:uniy}
If \(Y_{\delta}(\mbt)\) converges uniformly to \(Y(\mbt)\),
as defined in (\ref{eq:unifc}), then
\begin{eqnarray*}
  \lim_{\delta \to 0} \mP \left( \supt Y_{\delta}(\mbt) \ge
  c \right) = \mP \left( \supt Y(\mbt) \ge c \right) \quad \mbox{for any}
\quad c,
\end{eqnarray*}
where the right hand side is continuous.
\end{lemma}
{\em Proof.} For any \(\epsilon > 0\),
\begin{eqnarray*}
  \mP \left( \supt Y_{\delta}(\mbt) \ge c \right) \ge \mP \left(
  \supt  Y(\mbt) \ge c + \epsilon \right) - \mP \left( 
  \supt|Y_{\delta}(\mbt) - Y(\mbt)| > \epsilon \right).  
\end{eqnarray*}
Consequently, $\liminf_{\delta \to 0} \, \mP \left( 
\supt Y_{\delta}(\mbt) \ge c \right) \ge \mP \left( \supt Y(\mbt) \ge
c + \epsilon \right)$. However, since \(\epsilon\) is arbitrary, 
\begin{eqnarray*}
  \liminf_{\delta \to 0} \, \mP \left( \supt Y_{\delta}(\mbt)
  \ge c \right) \ge \mP \left( \supt Y(\mbt) \ge c \right).
\end{eqnarray*}
By a similar argument, it follows that
\begin{eqnarray*}
  \limsup_{\delta \to 0} \, \mP \left( \supt Y_{\delta}(\mbt)
  \right) \ge \mP \left( \supt Y(\mbt) \ge c \right)
\end{eqnarray*}
which completes the proof. \hfill \rule{2mm}{2mm}

\noindent {\em Proof of Theorem \ref{th:maxsstar}.}  First, convergence of
finite-dimensional distributions is a consequence of the multivariate
central limit theorem. Since a linear interpolant is always maximized
at one of the knots, this implies that the theorem holds for a linear
interpolant:
\[
  \lim_{n \to \infty} \mP \left( \supt
  S_{\delta}^{\star}(\mbt) \ge c \right) 
    = \mP \left( \supt Z_{\delta}(\mbt) \ge c \right)
\]
for any \(\delta > 0\). For any \(\epsilon > 0\), Lemma \ref{lem:unibd}
implies that
\begin{eqnarray*}
  \mP \left( \supt S^{\star}(\mbt) \ge c \right)
    &\le& \mP\left( \supt S_{\delta}^{\star}(\mbt) \ge c - \epsilon
    \right) + \mP \left( \supt | S^{\star}(\mbt) -
      S_{\delta}^{\star}(\mbt) | \ge \epsilon \right) \\
    &\le& \mP \left( \supt S_{\delta}^{\star}(\mbt) \ge c - \epsilon
    \right) + \frac{\delta}{\epsilon^2} \, \mE \int_{\mbT}
       \paren{ \frac{\partial}{\partial
           \mbt}S_{\delta}^{\star}(\mbt)}^2 \, d \mbt \\  
    &=& \mP \left( \supt S_{\delta}^{\star}(\mbt) \ge c - \epsilon
    \right) + \frac{\delta}{\epsilon^2} \, \mE \int_{\mbT}
       \left[ Z_{\delta}^{\prime}(\mbt) \right]^2 \, d \mbt,
\end{eqnarray*}
where $Z_{\delta}^{\prime}(\mbt) = \partial Z_{\delta}(\mbt)/\partial
\mbt$. The last equality follows from the fact that \(Z_{\delta}\) and
\(S_{\delta}^{\star}\) have the same covariance function. Assumption 
\ref{as:diffpsi} implies that the expectation is finite. From the
convergence of finite-dimensional distributions, it follows that 
\begin{eqnarray*}
  \limsup_{n \to \infty} \, \mP \left( \supt S^{\star}(\mbt) \ge
    c \right) \le \mP \left(\supt Z_{\delta}(\mbt) \ge c -
    \epsilon \right) + \frac{\delta}{\epsilon^2} \, \mE \int_{\mbT}
       \left[ Z_{\delta}^{\prime}(\mbt) \right]^2 \, d \mbt.
\end{eqnarray*}
First, let \(\delta \to 0\) and apply Lemma \ref{lem:uniy} to
\(Z_\delta\). Next, let \(\epsilon \to 0\) to obtain 
\begin{eqnarray*}
  \limsup_{n \to \infty} \, \mP \left( \supt S^{\star}(\mbt) \ge c
  \right) \le \mP \left(\supt Z(\mbt) \ge c \right).
\end{eqnarray*}
A similar argument shows that
\begin{eqnarray*}
  \liminf_{n \to \infty} \, \mP \left( \supt S^{\star}(\mbt) \ge
  c \right) \ge \mP \left( \supt Z(\mbt) \ge c
  \right)
\end{eqnarray*}
which completes the proof. \hfill \rule{2mm}{2mm}

\noindent {\em Proof of Theorem \ref{th:sig3d}.}  We assume the regularity
conditions 1 to 4 in \cite{adler:00}. The integral in
(\ref{eq:iprob}) can be expressed as
\begin{eqnarray}
   \int_{c^2}^{\infty} \mP\left(\supt \, 
    \ip{\bU_J}{\mbx_J(\mbt)} \ge \frac{c}{\sqrt{y}} \right)
    h_J(y) \, dy &=& \int_{c^2}^{\frac{c^2}{w_0}} \mP\left(\supt \, 
    \ip{\bU_J}{\mbx_J(\mbt)} \ge \frac{c}{\sqrt{y}} \right)
    h_J(y) \, dy \nonumber \\ 
   && + \int_{\frac{c^2}{w_0}}^{\infty} \mP\left(\supt \, 
    \ip{\bU_J}{\mbx_J(\mbt)} \ge \frac{c}{\sqrt{y}} \right)
    h_J(y) \, dy, \nonumber \\
    \label{eq:intsplit}
\end{eqnarray}
where $w_0 = (1 - r_0^2/2)$ and $r_0$ is the critical radius of
the tube. The volume-of-tube formula given in (\ref{eq:volexp}) is
exact when $y \in [c^2, c^2/w_0]$ and it is only approximate when $y
\in [c^2/w_0, \infty)$. In the former case, from (\ref{eq:tuberat}) 
\begin{eqnarray*}
   \mP\left(\supt \, \ip{\bU_J}{\mbx_J(\mbt)} \ge \frac{c}{\sqrt{y}}
   \right) = \sum_{t = 0}^d \frac{\zeta_t^J}{A_{d + 1 - t}} \mP
   \left[ B_{(d + 1 - t)/2, (J - d - 1 + t)/2} \ge w^2 \right].
\end{eqnarray*}
We express the first integral in (\ref{eq:intsplit}) as $F(c^2) -
F(c^2/w_0)$, where
\begin{eqnarray*}
  F(x) = \sum_{t = 0}^d \frac{\zeta_t^J}{A_{d + 1 - t}}
  \int_{x}^{\infty} \mP \left[ B_{(d + 1 - t)/2, (J - d - 1 + t)/2}
  \ge w^2 \right] h_J(y) \, dy.
\end{eqnarray*}
Note that the second integral in (\ref{eq:intsplit}) is $\geq 0$
providing a lower bound. Furthermore, 
\begin{eqnarray*}
  \int_{c^2/w_0}^{\infty} \mP\left(\supt \, 
    \ip{\bU_J}{\mbx_J(\mbt)} \ge \frac{c}{\sqrt{y}} \right)
    h_J(y) \, dy \leq \int_{c^2/w_0}^{\infty} h_J(y) dy
    = \mP\left(\chi^2_J \geq \frac{c^2}{w_0} \right).
\end{eqnarray*}
Therefore, $F(c^2) - F(c^2/w_0) \leq \mP\left[\supt \, Z_J(\mbt) \geq
c \right] \leq F(c^2) - F(c^2/w_0) + \mP\left(\chi^2_J \geq c^2/w_0
\right)$. As $c \to \infty$, $F(c^2) - F(c^2/w_0) \approx
F(c^2)$. Therefore, $\mP\left(\supt \, Z_J(\mbt) \geq
c \right) \approx F(c^2)$ as $c \to \infty$. By performing the
integration in $F(c^2)$, it follows that
\begin{eqnarray*}
  \mP\left(\supt \, Z_J(\mbt) \geq c \right) = \sum_{t = 0}^d
  \frac{\zeta_t^J}{A_{d + 1 - t}} \mP\left( \chi^2_{d + 1 - t} \ge
  c^2\right) + o[c^{-1} \exp(-c^2/2)] \quad \mbox{as} \quad c \to
  \infty. 
\end{eqnarray*}

When the Karhunen-Lo\`{e}ve expansion is infinite, the above result
for the truncated Gaussian random field $Z_J(\mbt)$ is extended by
letting $J \to \infty$ as follows. Uniform convergence of the
Karhunen-Lo\`{e}ve expansion implies that \(Z_J(\mbt) \longrightarrow
Z(\mbt)\) uniformly and hence
\begin{eqnarray}
  \label{eq:Zconv}
  \mP \left( \supt Z_J(\mbt) \ge c \right) \longrightarrow
  \mP \left( \supt Z(\mbt) \ge c \right)  \quad \mbox{as} \quad J
  \to \infty. 
\end{eqnarray}

The volume-of-tube formula given in (\ref{eq:volexp}) is in terms of
$\zeta^J_t$; however, as $J \to \infty$ and for $t = 0, \ldots, d$,
$\zeta^J_t \to \zeta_t$, the corresponding geometric term found via
$\rho(\mbt, \mbt)$, Therefore the result (\ref{eq:tubechis})
holds. For example, the expression for \(\k_0 \equiv \zeta_0\) is
derived by approximating \(\M\) by a series of short line segments to
obtain
\[
  \k_0 = \int_{\mbt} \Big| \det \left[\nabla_1 \nabla_2^T \rho(\mbt,
  \mbt) \right]\Big|^{1/2} d \mbt. \quad \quad \quad \quad 
\hfill \rule{2mm}{2mm}
\]
\remark We take sufficiently large $J$ so that the relation
(\ref{eq:Zconv}) holds. In practice, it is not necessary to employ a
truncated covariance function (\ref{eq:rhoJ}) that requires
specification of $J$ and the manifold $\M$. Our calculations are
carried out in terms of the covariance function $\mC(\mbt,
\mbt^{\dag})$. In effect, knowledge of $J$ and the specification of
$\M$ does not arise in practice. 

\setcounter{equation}{0}
\def\theequation{A.\arabic{equation}}

\newpage
%\vspace{0.2in}
{\bf Appendix A: Explicit Expressions for Geometric Constants 
in (\ref{eq:tubechis})}
%\vspace{0.05in}

We consider finite Karhunen-Lo\`{e}ve expansion with $J$ terms in
deriving the geometric constants. As a first step, we partition the
manifold $\M$, correspondingly the tube $\T(r, \M)$ and the parameter
space \(\mbT\) into various boundary regions. First, each point in
$\T(r, \M)$ is linked to a point in \(\M\) by a perpendicular
projection. Correspondingly, each point in \(\M\) is linked to a set
of points in $\T(r, \M)$. Second, partition \(\M\) into regions
\(\M_0, \ldots, \M_d\) based on the dimension of the linked sets,
where \(\M_0\) represents the main part of the manifold and \(\M_1,
\ldots, \M_d\) represent boundary regions. For example, when \(d =
1\), \(\M_1\) corresponds to the two end-points and \(\M_0\)
corresponds to the rest of the tube (see Fig.\
\ref{fig:tube}). If \(d = 2\), manifold \(\M\) is a polygon so that
\(\M_2\) represents the corners, \(\M_1\) the edges and \(\M_0\) the
interior. In effect, for a $d$-dimensional manifold $\M$, we can
partition both $\T(r, \M)$ and the space \(\mbT\) into \((d + 1)\)
regions to express $\vartheta(r, \M) = \V_0 + \V_1 + \cdots +
\V_d$. The main part of the tube can be represented as 
\begin{equation}
  \left[(1 + \|\mbtu\|^2)^{-1/2} \, \left(\mbx(\mbt) +
  \bQ(\mbt) \, \mbtu \right)\!: \mbt \in \mbT, \|\mbtu\| \le
  \tau_0 \right], \label{eq:reptu} 
\end{equation}
where $\tau_0 = \sqrt{1 - w^2}/w$, \(\bQ(\mbt)\) is an
orthonormal basis matrix for the normal space at
\(\mbx(\mbt)\). Provided that this transformation is one-to-one, the
volume \(\V_0\) can be expressed as $\V_0 = \int_{\mbt} \int_{\mbtu}
\Big|\det[ \bJ(\mbt, \mbtu)] \Big| \, d \mbt \, d \mbtu$, where
\(\bJ(\mbt, \mbtu)\) is the Jacobian of the representation
(\ref{eq:reptu}). The determinant of the Jacobian can be expressed as 
$\det[\bJ(\mbt, \mbtu)] = P_{\mbt}(\mbtu)(1 + \|\mbtu\|^2)^{-n/2}$, 
where \(P_{\mbt}(\mbtu)\) is a $d$th degree polynomial in \(\mbtu\)
with coefficients depending on \(\mbt\). This representation allows
the integral defining \(\V_0\) to be split into its \(\mbt\) and
\(\mbtu\) components, leading to a finite series expansion, for a
truncated $Z_J(\mbt)$, 
\[
  \V_0 = \sum_{t = 0}^d \k_t \frac{2 A_J}{A_{t + 1} A_{d + 1 - t}}
  \mP\left[B_{(d + 1 - t)/2, (J - d - 1 + t)/2} \ge
  w^2 \right], 
\]
where, \(\k_t\) are the polynomial coefficients integrated over
\(\M\) for even-order $t$ and the partial beta terms arise from 
integrating the \(\tau\) parts. Odd-order terms integrate to 0 by
symmetry; therefore, we set \(\kappa_t = 0\) when \(t\) is odd.  Recall
that \(A_t = 2 \pi^{t/2}/\Gamma(t/2)\) is the $(t - 1)$-dimensional
volume of the unit sphere $\S^{(t - 1)}$ in $\R^t$. The first constant
\(\k_0\) is the \(d\)-dimensional volume of the manifold $\M$,
represented in terms of the covariance function, expressed as 
\begin{equation}
   \k_0 = \int_{\mbT}
   \mC(\mbt, \mbt)^{-(d + 1)/2} \; \Bigg| \det \left[ \begin{matrix}
       \mC(\mbt, \mbt^{\dag}) & \nabla_2^T \, \mC(\mbt, \mbt^{\dag})
       \cr \nabla_1 \, \mC(\mbt, \mbt^{\dag}) & \nabla_1 \nabla_2^T \,
       \mC(\mbt, \mbt^{\dag}) \end{matrix}
	 \right]\Bigg|_{\mbt^{\dag} = \mbt}^{1/2} \, d \mbt,  
    \label{eq:multik0}
\end{equation}
where \(\nabla_1\) and \(\nabla_2\) denote vectors of partial
derivative operators with respect to the components of \(\mbt\) and
\(\mbt^{\dag}\) respectively. The geometric constant $\kappa_2$ is the
measure of curvature of $\M$.

The process for handling boundary corrections is similar.  To compute
the main boundary corrections, represent the half-tubes around
boundaries in a form similar to (\ref{eq:reptu}), with
\(\bQ(\mbt)\) supplemented by a vector tangent to $\M$ but
normal to $\partial \M$, the boundary of $\M$. The vector \(\mbtu\) is
then restricted to a half-sphere. Following the derivation of
\cite{weyl:39}, we obtain a series of the form, for truncated
$Z_J(\mbt)$,    
\begin{eqnarray*}
  \V_1 = \sum_{t = 0}^{d - 1} \ell_t \, \frac{A_J}{A_{t + 1}
     A_{d - t}} \mP \left[ B_{(d - t)/2, (J - d +
     t)/2} \ge w^2 \right],
\end{eqnarray*}
where \(\ell_t\) terms are the integrals of polynomial
coefficients. The first term, \(\ell_0\) is the \((d -
1)\)-dimensional volume of $\partial \M$ which has a form similar to
(\ref{eq:multik0}), summed over each of the boundary faces. It is
important to note that odd order terms no longer disappear; $\ell_1$
is a measure of rotation of $\partial \M$ and $\ell_2$ is a measure of
curvature similar to $\kappa_2$. Similarly, at corners where two
boundary faces meet, we can represent
\begin{eqnarray*}
  \V_2 = \sum_{t = 0}^{d - 2} \nu_t \frac{A_J}{A_{t + 1} A_{d - 1 -
  t}} \, \mP\left[B_{(d - 1 - t)/2, (J - d + 1 +
  t)/2} \ge w^2 \right], 
\end{eqnarray*}
where \(\nu_0\) measures the rotation angles in the regions of
$\partial^2 \M$ (the boundary of $\partial \M$) where two boundary
faces meet and $\nu_1$ is a combination of rotation angles and
rotation of the edges. Currently, our software library enables
computing all the terms given in (\ref{eq:tubechis}); effectively
yielding a complete implementation of the tube formula up to \(d =
3\). To the best of our knowledge, there exist no method for general
implementation of higher-order terms with boundary corrections.

\remark When \(d = 2\), the fourth order coefficients are \(\ell_2 =
\nu_1 = m_0 = 0\). Additionally, the Euler-Poincare characteristic
\citep{knowles:89} satisfies $\kappa_2 + \ell_1 + \nu_0 = 2 \pi
{\mathcal E} - \kappa_0$ eliminating the need to compute \(\kappa_2\),
\(\ell_1\) and \(\nu_0\) directly. The Euler-Poincare characteristic
is the number of pieces making up the manifold, minus the number of
holes. When \(\mbT\) is a compact as well as a convex set and
\(\mC(\mbt, \mbt) > 0\) for all \(\mbt\) then \({\mathcal E} = 1\). 

Combining the above results together, the tube formula, up to fourth
order terms, can be expressed as 
\begin{eqnarray}
  \vartheta(r, \M) &\approx& \V_0 + \V_1 + \V_2 + \V_3 \nonumber \\
  &=& \k_0 \frac{A_J}{A_{d + 1}} \mP \left[ B_{(d + 1)/2,
  (J - d - 1)/2} \ge w^2 \right] + \ell_0 \frac{A_J}{2 A_d}
  \mP\left[ B_{d/2, (J - d)/2} \ge w^2 \right]
  \nonumber \\ && \; + \left(\k_2 + \ell_1 + \nu_0 \right)
  \frac{A_J}{2 \pi A_{d - 1}} \mP \left[ B_{(d - 1)/2,
  (J - d - 1)/2} \ge w^2 \right] \nonumber \\ && \; +
  \left(\ell_2 + \nu_1 + m_0 \right) \frac{A_J}{4 \pi A_{d - 2}}
  \mP\left[ B_{(d - 2)/2, (J - d - 2)/2} \ge w^2
  \right], \label{eq:vol} 
\end{eqnarray}
where \(m_0\) measures the size of wedges at corners where three
boundary faces of \(\M\) meet. After completing evaluation of all
terms leads to a series, 
\begin{eqnarray}
  \label{eq:volexp}
  \vartheta(r, \M) \approx \sum_{t = 0}^d \zeta_t^J \frac{A_J}{A_{d
  + 1 - t}} \mP \left[ B_{(d + 1 - t)/2, (J - d - 1 + t)/2} \ge w^2
  \right].
\end{eqnarray}
The dominant term \(\zeta_0^J\) can be expressed as 
\begin{eqnarray}
  \label{eq:zeta}
  \zeta_0^J = \int_{\mbT} \Big\| \frac{\partial}{\partial \mbt}
  \mbx(\mbt) \Big\| \, d \mbt = \int_{\mbt} \Big| \det
  \left[\nabla_1 \nabla_2^T \rho_{J}(\mbt, \mbt) \right]\Big|^{1/2} d
  \mbt, 
\end{eqnarray}
where \(\nabla_1\) and \(\nabla_2\) are partial derivative operators
with respect to the first and second arguments of \(\rho_J(\cdot,
\cdot)\), respectively. 

The following correspondence (up to \(t = 3\)) holds: $\zeta_0 = \k_0,
\zeta_1 = \ell_0/2, \zeta_2 = (\k_2 + \ell_1 + \nu_0)/(2 \pi)$ and
$\zeta_3 = (\ell_2 + \nu_1 + m_0)/(4 \pi)$. The tube formula is  exact
for tubes with radius $r \leq r_0$, the critical radius. 
\bibliographystyle{apalike}
\bibliography{tube-mixv2}

\end{document}